\documentclass{article}
\usepackage[utf8]{inputenc}
\usepackage[english]{babel}
\usepackage{lmodern}
\usepackage[babel=true]{microtype}
\usepackage{amsmath,amssymb}
\usepackage{amsfonts}
\usepackage{amsthm}
\usepackage{vmargin}
\usepackage{csquotes}
\usepackage{bbm}
\usepackage{extarrows}
\usepackage{titling}
\usepackage{mathtools} 
\usepackage[style=alphabetic, backend=biber,]{biblatex}
\DeclareLabelalphaNameTemplate{
  \namepart[compound=true]{family}
}

\usepackage{titlesec}

\titleformat{\section}[block]
  {\normalfont\Large\bfseries\centering}
  {\Roman{section}.}{0.5em}{}

\titleformat{\subsection}[block]
  {\normalfont\large\bfseries\centering}
  {\Roman{section}.\arabic{subsection}.}{0.5em}{}

\usepackage{pgf,tikz,tikz-cd}
\usetikzlibrary{arrows}
\usetikzlibrary{arrows.meta}

\usepackage{hyperref}
\hypersetup{
	colorlinks=true,
	linkcolor=blue,
	filecolor=blue,
	citecolor = blue,      
	urlcolor=blue,
}

\theoremstyle{plain}
\newtheorem{theorem}{Theorem}[subsection]
\newtheorem{definition}[theorem]{Definition} 

\newtheorem{proposition}[theorem]{Proposition}

\newtheorem{corollary}[theorem]{Corollary}
\newtheorem{lemma}[theorem]{Lemma}

\newtheorem{conjecture}[theorem]{Conjecture}
\newtheorem{question}[theorem]{Question}

\theoremstyle{remark}
\newtheorem{remark}[theorem]{Remark}

\newtheorem{notation}[theorem]{Notation}

\newcommand{\Hom}{\mathrm{Hom}}

\newcommand{\Spec}{\mathrm{Spec}\,}

\newcommand{\CH}{\mathrm{CH}}
\newcommand{\id}{\mathrm{id}}
\newcommand{\Corr}{\mathrm{Corr}}
\newcommand{\longhookrightarrow}{\lhook\joinrel\longrightarrow}

\pretitle{\begin{center}\large\bfseries}
\posttitle{\end{center}}

\preauthor{\begin{center}\large}
\postauthor{\end{center}\vspace{-5em}}

\title{\MakeUppercase{Homogeneous spaces over an abelian variety}}
\author{\textsc{Margot Bruneaux}\thanks{The author was supported by the project "Group schemes, root systems, and related representations" founded by the European Union - NextGenerationEU through Romania’s National Recovery
and Resilience Plan (PNRR) call no. PNRR-III-C9-2023- I8, Project CF159/31.07.2023, and coordinated by the Ministry of Research, Innovation and Digitalization (MCID) of Romania, and by the LABEX MILYON (ANR-10-LABX-0070) of Université de Lyon, within the program "France 2030" (ANR-11-IDEX-0007) operated by the French National Research Agency (ANR).}}
\date{}

\setmarginsrb{3cm}{3cm}{3cm}{2cm}{0cm}{0cm}{0cm}{1cm}

\begin{document}

\maketitle
\begin{abstract}
In this paper, we study a question of Colliot-Thélène and Iyer concerning the existence of rational sections in families of homogeneous spaces over an abelian variety, after base change by a suitable étale isogeny of the abelian variety. Assuming characteristic zero and that the homogeneous spaces arise from connected reductive groups, the problem is reformulated in terms of torsors under reductive groups over an abelian variety $A$.

Building on work of Moonen and Polishchuk, we construct a filtration on the motive of a Jacobian variety to analyze the action of isogenies on unramified cohomology and Witt groups. This approach allows for a positive response to the question for reductive groups whose root data do not contain a factor of type~$E_8$ when $\dim A > 2$ and $\mathrm{cd}(k) \leqslant 1$, and for all reductive groups when $\dim A = 2$ and $k$ is algebraically closed.

\bigskip

\noindent\textsc{2020 Mathematics Subject Classification:} 14L15, 14K05, 14K02, 14C15, 14F42, 11E81

\smallskip

\noindent \textsc{Keywords:} torsors, homogeneous spaces, reductive groups, abelian varieties.

\end{abstract}

\makeatletter
\@addtoreset{theorem}{section}
\renewcommand{\thetheorem}{\thesection.\arabic{theorem}}
\makeatother

\section*{\normalfont \textsc{Introduction}}

The study of families of projective homogeneous spaces over abelian varieties is motivated by questions of potential density and the geometry of varieties with nef tangent bundle. Colliot-Thélène and Iyer proved that, over a number field, potential density on the base of a smooth projective family of homogeneous spaces implies potential density on its total space. For an abelian base, this gives potential density after a finite extension of the ground field.

Seeking a more geometric explanation for their result, and motivated by the Campana–Peternell conjecture, they asked whether, over an algebraically closed field, the obstruction to a rational section of such a family can always be eliminated by passing to a finite étale cover of the abelian variety:

\begin{question}{\cite[Question~3.4]{ColliotThélèneIyer2011}}
Let $A$ be an abelian variety over an algebraically closed field~$k$.
Let $X \longrightarrow A$ be a smooth, projective family of homogeneous spaces
of connected linear algebraic groups.
Does there exist a finite étale map $B \longrightarrow A$, with $B$ connected,
such that the base change $X \times_A B \longrightarrow B$ admits a rational section?
\end{question}

In this paper, we give a nearly complete affirmative answer to this question
in characteristic~$0$, working moreover over base fields of cohomological dimension at most~$1$.
Since the unipotent radical of an algebraic group acts trivially on a projective homogeneous space,
we may assume without loss of generality that the homogeneous space arises from a connected reductive group.
Our main result gives a positive answer whenever the root data of the groups acting on the fibers
contain no factors of type~$E_8$:

\begin{theorem}{[Theorem \ref{Theorem_homogeneous}]}
Let $A$ be an abelian variety over a field~$k$ of cohomological dimension at most~$1$,
of characteristic~$0$, and containing a square root of~$-1$.
Let $U \subset A$ be an open subset that contains all points of codimension~$1$, with $0 \in U$.
Let $X \longrightarrow U$ be a smooth, projective family of homogeneous spaces
under connected reductive algebraic groups.
Assume that the type of $X$ has no factor of type $E_8$. Then there exists an étale isogeny $f : A \to A$ such that the pullback family
\[
f^*X \to f^{-1}(U)
\]
admits a rational section.
\end{theorem}

\medskip

In Section~\ref{TheconjectureforBorelvarieties},
we use a theorem of Demazure~\cite[Prop.~4]{Demazure1977}
to reformulate the problem in the natural language of Galois cohomology,
as a question about torsors under reductive groups:

\begin{question} \label{conjecture_torsor_intro}
Let $A$ be an abelian variety over a field~$k$, and let $U \subset A$ be an open subset
that contains all points of codimension~$1$, with $0 \in U$.
Let $G$ be a reductive group defined over~$U$,
and let $X \longrightarrow U$ be a $G$-torsor trivial at the origin.
Does there exist an étale isogeny $f : A \longrightarrow A$ such that the pullback
$f^*X \longrightarrow f^{-1}(U)$ is Zariski-locally trivial?
\end{question}

Berkovich~\cite[p.~182]{Berkovich1973} gave an affirmative answer in the case
$G = \mathrm{PGL}_n$, $k$ algebraically closed, and $U = A$,
using cohomological invariants from the Brauer group together with the structure of
$H^i(A, \mathbb{Z}/n)$, which is fully determined by $H^1(A, \mathbb{Z}/n)$
(see Proposition~\ref{etale_finite_cohomology}).
We extend the use of cohomological invariants, a theory developed by Serre in~\cite{Garibaldi2003},
to give a positive answer to Question~\ref{conjecture_torsor_intro}
when the fiber $G_0$ at the origin is quasi-split with root datum containing no $E_8$ factor,
and when $k$ contains a square root of $-1$
(see Theorem~\ref{Theorem_torsor}).
By Steinberg's theorem~\cite[Cor.~5.2.6]{Gille2019},
the quasi-split assumption is automatic when $\mathrm{cd}(k) \leqslant 1$,
and in particular when $k$ is the function field of an algebraic curve
over an algebraically closed field.

\medskip

The exclusion of the $E_8$ case stems from a limitation of the method:
cohomological invariants are effective when there exist sufficiently many of them
to detect the triviality of torsors under almost simple quasi-split groups,
which is the case for all types except $E_8$.
For groups of type $E_8$, no such invariants are currently known,
and the question remains open.
For abelian surfaces over an algebraically closed field,
we circumvent this difficulty by appealing to a theorem
of de~Jong, He, and Starr~\cite[Theorem~1.4]{Jong2011},
which yields a positive answer regardless of the root datum.

\medskip

The cases of split orthogonal groups of types $B_n$ and $D_n$ are of particular
importance for the resolution of this problem.
To handle the question for this class of groups, we make use of the Milnor invariants
via the \emph{Pardon spectral sequence} (Section~\ref{Pardon_spectral_sequence}).
A key difficulty is that Milnor invariants are not globally defined on the variety,
but are only available in the unramified setting,
which forces one to work systematically with unramified cohomology throughout.
This leads to the following result on the behavior of unramified cohomology under isogenies:

\begin{proposition}[Proposition \ref{cancelling_unramified_coho}]
Let $k$ be a field of characteristic zero, let $A$ be an abelian variety over $k$,
and let $n \in \mathbb{N}^*$.
Then there exists an étale isogeny $f$ of $A$ such that,
for every $\alpha \in H^i_{\mathrm{nr}}\!\bigl(A, \mathbb{Z}/n\bigr)$,
the pullback $f^*(\alpha)$ is constant,
that is, it comes from an element of $H^i_{\mathrm{nr}}\!\bigl(k, \mathbb{Z}/n\bigr)$.
\end{proposition}

This result, which allows one to trivialize all unramified cohomology classes
after a suitable isogeny, is of independent interest
and plays a central role in the treatment of orthogonal groups.
We also study the action of the multiplication-by-$n$ morphism~$[n]$
on the unramified Witt group of an abelian variety
(see Section~\ref{unramified_witt_group}).

\medskip

To establish these results,
we require a finer understanding of the Voevodsky motive of an abelian variety
with integral coefficients.
While the rational structure is well understood
via the work of Ancona--Huber--Pépin Lehalleur~\cite{AnconaHuberPepinLehalleur2016},
the integral situation is considerably more subtle.
In the case where $J$ denotes the Jacobian of a smooth, projective, connected curve over~$k$,
and $S$ is a smooth quasi-projective $k$-scheme,
Moonen and Polishchuk~\cite{MoonenPolishchuk2010} provide a detailed analysis
of the Chow groups~$\CH_*(J \times_k S)$.
Building on their work,
we construct a filtration of the dual cohomological motive of a Jacobian
with favorable compatibility properties under multiplication maps~$[n]$
(see Proposition~\ref{defining_filtration} and Corollary~\ref{f_p_maps}).
Via Voevodsky's functor~\cite[Prop.~2.1.4]{Voevodsky2000},
we obtain the following structural result (see Proposition~\ref{action_n_DM}):

\smallskip
\emph{There exists a decomposition $M(J) = M(k) \oplus M'$ in $\mathsf{DM_{gm}}(k)$
such that $M'$ is stable under $\mathsf{M_{gm}}([n])$,
and the restriction of $\mathsf{M_{gm}}([n])^{3\dim J}$ to $M'$ is of the form $nf$,
where $f$ is an endomorphism of~$M'$.}


\medskip

Finally, the \emph{Campana--Peternell conjecture} predicts that,
up to a finite étale cover,
every smooth complex projective variety with a nef tangent bundle
is a homogeneous space over an abelian variety, more precisely, that a Fano variety with nef tangent bundle
is a projective homogeneous space of a linear algebraic group
(over an algebraically closed field of characteristic zero),
a statement known in dimension less than~$4$.
The existence of rational sections for smooth homogeneous spaces over abelian varieties
is closely related to the geometry of such varieties,
and in Section~\ref{Generalisation}
we extend our results to smooth projective homogeneous spaces over abelian varieties,
contributing to this conjecture.

\medskip

\noindent \textbf{Acknowledgements.} I would like to sincerely thank Philippe Gille for his guidance and support, without which this article would not have been written. I am also grateful to Paolo Bravi for our discussions and for sharing his perspectives on the problem. I would like to thank Federico Scavia, Swann Tubach, Giuseppe Ancona, Baptiste Calmès, János Kollár, and Simone Diverio for helpful discussions. Finally, I thank Jean-Louis Colliot-Thélène for pointing out that the hypothesis on bounded cohomological dimension of the base field in Theorem~\ref{Theorem_torsor} can be removed.

\makeatletter
\@addtoreset{theorem}{subsection}
\renewcommand{\thetheorem}{\thesubsection.\arabic{theorem}}
\makeatother

\subsection*{\normalfont \textsc{Notation}}

Let $k$ be a field, let $\mathsf{SmProj}(k)$ denotes the category of smooth projective schemes over $k$ and let $\mathsf{Sm}(k)$ denotes the category of smooth schemes over $k$.

Let $\mathsf{CSmProj}(k)$ denotes the category of correspondences of smooth projective schemes over $k$ (see \cite[§4.1]{André2004}).

Let $\mathcal{M}(k)$ be the category of pure  Chow motives over $k$ with integral coefficients (see \cite[§4.1]{André2004}).

Let $\mathcal{M}^{\mathrm{eff}}(k)$ be the category of effective Chow motives over $k$ with integral coefficients (see \cite[§4.1]{André2004}).

Let $\mathsf{DM}(k)$ denote Voevodsky's category of motives over $k$ and $\mathsf{DM_{gm}}(k)$ denote Voevodsky's category of geometric motives (see \cite{Voevodsky2000}).

Let $\mathsf{GrAb}$ be the category whose objects are graded $\mathbb{Z}$-modules and whose morphisms are morphisms of $\mathbb{Z}$-modules (not necessarily graded).  

If $X \in \mathsf{SmProj}(k)$, we denote by $h(X)$ its Chow motive, by $h(X)^{\vee}$ its dual motive, and by $M(X)$ its motive in $\mathsf{DM}_{\mathrm{gm}}(k)$. Recall that for an equidimensional smooth projective scheme $X$ of dimension $d$, we have
\[
h(X)^{\vee} = h(X)(d).
\]
For more details on duals in $\mathcal{M}(k)$, we refer the reader to \cite[§4.1]{André2004}.

We have a natural map $\alpha$ from the big Zariski site to the big étale site of a scheme $X$. We denote by $(\alpha^*_X,\alpha_{X,*})$ the associated maps of topoi. When $X=\Spec k$, we simply write $(\alpha^*,\alpha_{_*})$.

For an étale sheaf $\mathcal{F}$ on the big étale site of $\Spec k$, we denote by $\mathcal{H}^i(\mathcal{F})$ the sheaf $\mathcal{R}^i\alpha_* \mathcal{F}$.

For a constructible sheaf $\mathcal{F}$  of $n$-torsion on the small étale site of a scheme $X$, we put $\mathcal{F}(i)= \mathcal{F}\otimes\mathbb{Z}/n (i)$ for every $i \in \mathbb{Z}$.

Let $R$ be a discrete valuation ring (DVR) containing a field of characteristic different from $2$, and let $K$ be a field of characteristic different from $2$. We denote by $W(R)$ and $W(K)$ the Witt groups of quadratic forms over $R$ and $K$, respectively (see \cite{Merkurjev2008} for the definitions).

We denote by $I(K) \subset W(K)$ the \emph{fundamental ideal} of $W(K)$, that is, the kernel of the rank homomorphism
\[
\operatorname{rk} : W(K) \longrightarrow \mathbb{Z}/2\mathbb{Z},
\]
which assigns to each quadratic form its rank modulo $2$. Equivalently, $I(K)$ is the ideal of $W(K)$ generated by the classes of quadratic forms of even rank. For each integer $i \geqslant 0$, we denote by $I^i(K)$ the $i$-th power of this ideal.

Unless otherwise specified, all cohomology groups $H^i$ will be understood as étale cohomology groups $H_{\mathrm{\acute{e}t}}^i$, and all isogenies of abelian varieties will be taken to be morphisms of groups.

\section{A filtration on the dual motive of a Jacobian variety} 

In this section, we deduce a decomposition of the Chow motive of the Jacobian of a curve from a result of Moonen and Polishchuk on its Chow groups.

\begin{notation} \label{Notation_1}Let $k$ be a field.

\noindent Let $S$ be a smooth, projective $k$-scheme of dimension $d$.

\noindent Let $C$ be a smooth, projective, connected curve over $k$, equipped with a distinguished $k$-point $p_0$. 

\noindent For every $m \in \mathbb{N}$, denote by $C^{[m]}$ the $m$th symmetric power of $C$, with the convention that $C^{[m]} = \emptyset$ if $m < 0$.  

\noindent Let $J$ be the Jacobian variety of $C$ and $g$ its dimension.  

\noindent The $n$th symmetric power $C^{[m]}$ is a smooth $k$-scheme (see \cite[Prop.~3.2]{Milne1986}), and there is a natural morphism depending on $p_0$
\[
\sigma_m \colon C^{[m]} \longrightarrow J,
\]
called the Albanese map (see \cite[§5]{Milne1986}). The Albanese map from $C^{[m]}$ to $J$ induces a morphism $({\sigma_m})_*$  from $\CH_*(C^{[m]}\times_k S)$ to $\CH_*(J\times_k S)$. The sum of these morphisms induces a map from $\CH_*(C^{[\bullet]}/S)$ to $\CH_*(J\times_k S)$, which we shall denote by $(\widetilde{\sigma}_{S})_*$.

\noindent For every $n\in\mathbb{N}$, let $[n]: J \longrightarrow J$ denotes the multiplication by $n$. Since it is a proper map, it induces a $0$-graded endomorphism of the Chow group $\CH_*(J\times _k S)$ denoted by $[n]_{S,*}$.
    \end{notation}

In the setting Notation \ref{Notation_1}, our aim is to define a filtration on the dual  $h(J)^{\vee}$ of the Chow motive $h(J)$ of the Jacobian $J$. To achieve this, we use the filtration introduced by Moonen and Polishchuk in \cite{MoonenPolishchuk2010} on the Chow groups $\mathrm{CH}_*(J \times_k T)$ for every smooth and quasi-projective $k$-scheme $T$.
 This filtration is interesting because of its nice properties related to the action of $[n]_{*}$ on $h(J)^{\vee}$.
We first recall some results from \cite{MoonenPolishchuk2010}, then deduce the filtration and describe the action of $[n]_*$ on it using the Manin Identity principle.

\subsection{The inclusion of \texorpdfstring{$\CH_*(J\times_k S)$}{} in \texorpdfstring{$\CH_*(C^{[\bullet]}/S)$}{} and the induced filtration} \label{Moonen_Polishchuk_result}

One of the important results of \cite{MoonenPolishchuk2010} is the construction of an embedding
\[
\mathrm{CH}_*(J \times_k S) \longhookrightarrow \mathrm{CH}_*(C^{[\bullet]} \times_k S)
= \bigoplus_{m \in \mathbb{N}} \mathrm{CH}_*(C^{[m]} \times_k S).
\]
This inclusion also respects the algebra structures on $\CH_*(J \times_k S)$ and $\CH_*(C^{[\bullet]} \times_k S)$ induced by the Pontryagin product (see \cite[p.~805]{MoonenPolishchuk2010} for a definition of the Pontryagin product). However, the multiplicative structure of these Chow groups will not be relevant for us, so we will not provide further details on it in the following. The image of this inclusion is generated by sums of homogeneous elements in the image with respect to the $C^{[\bullet]}$-grading on $\CH_*(C^{[\bullet]} \times_k S)$, where $\CH_*(C^{[n]} \times_k S)$ is placed in degree $n$. Consequently, it induces a grading on $\CH_*(J\times_k S)$ which behaves nicely under the operator $[n]_{S,*}$.
We will only provide an outline of these results, established in \cite{MoonenPolishchuk2010}; for further details, we refer the reader to the original article. To obtain the embedding $\mathrm{CH}_*(J \times_k S) \longhookrightarrow  \mathrm{CH}_*(C^{[\bullet]} \times_k S)$, the idea is to construct a section of the map $(\widetilde{\sigma}_{S})_*$.

\subsubsection*{Operators on the Chow homology}

In his article \cite{Polishchuk2007}, the author defines a family of homogeneous operators on $\CH_*(C^{[\bullet]}\times_k S)$ for every smooth projective scheme $S$ over $k$. For every $i,j\in\mathbb{N}$ and $a\in\CH_*(C\times_k S)$, he defines an operator denoted $P_{i,j}(a):\CH_*(C^{[\bullet]}\times_k S)\longrightarrow\CH_*(C^{[\bullet+i-j]}\times_k S)$ of degree $i-j$. 
He also defines the divided powers of the operators $P_{0,n}(C\times_k S)$ and $P_{n,0}(C\times_k S)$. Instead of giving the general definition of these operators, we present only the expressions relevant to our discussion:

    \begin{align}\label{Operators}
\begin{split}
&\forall m\in \mathbb{Z},\ \forall x\in \CH_*(C^{[m+1]}\times_k S),\\
&P_{0,1}([p_0]_S)\big|_{\CH_*(C^{[m+1]}\times_k S)}(x)
= (i_m)_S^*(x)
\qquad \text{\cite[\S 3, p.~561]{Polishchuk2007}},\\[0.5em]
&\forall m\in \mathbb{Z},\ \forall j\in \mathbb{N},\ \forall x\in \CH_*(C^{[m]}\times_k S),\\
&P_{0,1}(C\times_k S)^{[j]}\big|_{\CH_*(C^{[m]}\times_k S)}(x)
= (\mathrm{pr}_2)_*\bigl(p_1^*([C^{[j]}]) \cdot \alpha_{j,m-j}^*(x)\bigr) \\
&\phantom{P_{0,1}(C\times_k S)^{[j]}\big|_{\CH_*(C^{[m]}\times_k S)}(x)}
= (\mathrm{pr}_2)_*\bigl(\alpha_{j,m-j}^*(x)\bigr)
\qquad \text{\cite[p.~551]{Polishchuk2007}}.
\end{split}
\end{align}

\noindent
where $i_m : C^{[m]} \to C^{[m+1]}$ denotes the natural inclusion (with respect to the first factor), 
$\alpha_{j,m-j} : C^{[j]}\times C^{[m-j]} \to C^{[m]}$ is the natural morphism, 
and $P_{0,1}(C\times_k S)^{[j]}$ denotes the $j$-th divided power of the operator $P_{0,1}(C\times_k S)$.

\medskip

\subsubsection*{Construction of the section}

On $\CH_*(C^{[\bullet]}/S)$, there is a natural grading defined as follows.

\begin{definition}\label{grading_def}
The \emph{$C^{[\bullet]}$-grading} on $\CH_*(C^{[\bullet]} \times_k S)$ is the grading for which $\CH_*(C^{[m]}\times_k S)$ is placed in degree $m$. 
\end{definition}

\begin{lemma}\label{K_S}
The submodule of $\CH_*(C^{[\bullet]}/S)$
\[
\mathbb{K}_S := \mathrm{Ker}\bigl(P_{0,1}([p_0]_S)\bigr) \cap \bigcap_{n \geq 1} \mathrm{Ker}\bigl(P_{0,1}(C\times_k S)^{[n]}\bigr)
\subset \CH_*(C^{[\bullet]}/S)
\]
is generated by its \emph{homogeneous elements} with respect to the $C^{[\bullet]}$-grading for which $\CH_*(C^{[m]}\times_k S)$ is in degree $m$.
\end{lemma}

\begin{proof}
Each of the operators $P_{0,1}([p_0]_S)$ and $P_{0,1}(C\times_k S)^{[n]}$ is homogeneous with respect to this grading. 
More precisely, $P_{0,1}([p_0]_S)$ is homogeneous of degree $-1$, i.e.\ it sends $\CH_*(C^{[m]}\times_k S)$ to $\CH_*(C^{[m-1]}\times_k S)$ for every $m\in \mathbb{Z}$. 
Similarly, for every $n\in \mathbb{N}_{>0}$, the operator $P_{0,1}(C\times_k S)^{[n]}$ is homogeneous of degree $-n$, i.e.\ it sends $\CH_*(C^{[m]}\times_k S)$ to $\CH_*(C^{[m-n]}\times_k S)$ for every $m\in \mathbb{Z}$. This proves the result.
\end{proof}

In \cite[Thm 3.4]{MoonenPolishchuk2010}, the following property is proven:

\begin{proposition} \label{K_S_CH_J}Let $k$ be a field, let $C$, $J$, $S$ be as in Notation \ref{Notation_1} and let $\mathbb{K}_S$ be as in Lemma \ref{K_S}.
    The homomorphism $({\widetilde{\sigma}}_{S})_* : \CH_*(C^{[\bullet]}/S) \longrightarrow \CH_*(J \times_k S)$ restricts to an isomorphism $\mathbb{K}_S \xlongrightarrow[]{\sim} \CH_*(J\times_k S)$. 
\end{proposition}
\begin{proof}
 This is the first part of Theorem 3.4 in \cite{MoonenPolishchuk2010}
\end{proof}

Let $S$ be a projective smooth $k$-scheme. The composition of the inverse isomorphism
$\mathbb{K}_S \xlongrightarrow{\sim} \CH_*(J\times_k S)$
with the inclusion
$\mathbb{K}_S \hookrightarrow \CH_*(C^{[\bullet]} \times_k S)$ factors through a map
\begin{equation}\label{s_S_def}
\widetilde{s}_S : \CH_*(J\times_k S) \longrightarrow \bigoplus_{m=0}^{2 \dim J + \dim S} \CH_*(C^{[m]} \times_k S)
\end{equation}
since, by \cite[Proposition 3.11 (i)]{MoonenPolishchuk2010}, $\mathbb{K}_S$ is contained in $\bigoplus_{m=0}^{2 \dim J + \dim S} \CH_*(C^{[m]} \times_k S)$.

\begin{proposition}\label{s_S}
Let $k$ be a field, and let $C$, $J$, and $S$ be as in Notation \ref{Notation_1}. The map $\widetilde{s}_S$ of \eqref{s_S_def} is an injective morphism from $\CH_*(J\times_k S)$ to $\CH_*(C^{[\bullet]} \times_k S)$ that preserves dimensions, i.e. it sends $\CH_n(J\times_k S)$ into $\CH_n(C^{[\bullet]} \times_k S)$ for every $n \in \mathbb{N}$. Moreover, its image is generated by homogeneous elements with respect to the grading on $\CH_*(C^{[\bullet]} \times_k S)$ defined in Definition~\ref{grading_def}.
\end{proposition}

\begin{proof}
As the composition of the inverse isomorphism
$\mathbb{K}_S \xlongrightarrow{\sim} \CH_*(J\times_k S)$
with the inclusion
$\mathbb{K}_S \hookrightarrow \CH_*(C^{[\bullet]} \times_k S)$ is injective, so is $\widetilde{s}_S$.

Since $\mathbb{K}_S$ is generated by its homogeneous elements with respect to the $C^{[\bullet]}$-grading on $\CH_*(C^{[\bullet]} \times_k S)$ by Lemma \ref{K_S}, the image of $\widetilde{s}_S$ is generated by homogeneous elements.

Moreover, since $(\widetilde{\sigma}_S)_*$ sends $\CH_n(C^{[\bullet]} \times_k S)$ into $\CH_n(J\times_k S)$ for every $n$, it follows that $\widetilde{s}_S$ sends $\CH_n(J\times_k S)$ into $\CH_n(C^{[\bullet]} \times_k S)$ for every $n$.
\end{proof}

\subsubsection*{A filtration of $\CH_*(J\times_k S)$}

This section $\widetilde{s}_S$ induces a $C^{[\bullet]}$-grading on $\CH_*(J \times_k S)$, and its image is generated by homogeneous elements with respect to the $C^{[\bullet]}$-grading on $\CH_*(C^{[\bullet]}/S)$. In \cite{MoonenPolishchuk2010}, the authors show that this grading is bounded and induces a filtration with nice properties related to the action of $[n]^*$.
\begin{samepage}
\begin{theorem}{\cite[Thm.~4]{MoonenPolishchuk2010}} \label{Theorem_4} Let $C$, $J$, $S$ be as in Notation \ref{Notation_1} and let $\widetilde{s}_S$ be as in (\ref{s_S_def}).  
\begin{enumerate}
    \item[\textup{(i)}]We have a decomposition
\[ \CH_*(J\times_k S) = \bigoplus_{m=0}^{2\dim J+ \dim S} \CH_*^{[m]}(J \times_k S) \]
where $x \in \CH_*^{[m]}(J\times_k S)$ if and only if ${\widetilde{s}}_S(\alpha) \in \CH_*(C^{[m]}\times_k S) \subset \CH_*(C^{[\bullet]}/S)$.
\item[\textup{(ii)}] For every $n \in \mathbb{Z}$, the associated descending filtration $\mathrm{Fil}^{\bullet}$ on $\CH_*(J \times_k S)$ is stable under the operators $[n]_{S,*}$ and $[n]_{S,*}$ acts on ${\mathrm{gr}}_{\mathrm{Fil}}^m$ as the multiplication by $n^m$.
\end{enumerate} 
\end{theorem}
\end{samepage}

\subsection{Motivic interpretation} \label{Motivic_interpretation}

 In this section, we present some consequences of the previous one. We will show that the filtration on $\CH_*(J \times_k S)$ is functorial in $S$, where morphisms between schemes are given by correspondences. This functoriality will allow us to lift results about Chow groups to the level of motives via the Manin Identity principle.

\subsubsection*{Manin Identity Principle}

Let us recall that we denote by $\mathsf{GrAb}$ the category whose objects are graded $\mathbb{Z}$-modules and whose morphisms are $\mathbb{Z}$-linear maps (not necessarily homogeneous).

\begin{definition}\label{def_omega} 
Let $k$ be a field. To every motive $M \in \mathcal{M}(k)$, we can associate a functor:
\[
\begin{tabular}{cccc}
$\omega_{M}$  : & $\mathsf{CSmProj}(k)^{\mathrm{op}}$ &  $\longrightarrow$  & $\mathsf{GrAb}$ \\
 & $S$ & $\longmapsto$ & $\displaystyle \bigoplus_{r \in \mathbb{Z}} \mathrm{Hom}_{\mathcal{M}(k)}(h(S),M(r))$.
\end{tabular}
\]
Let $r' \in \mathbb{Z}$,
$c \in \Corr^{r'}(S',S) := \CH^{\dim S' + r'}(S' \times_k S) \subset \Hom_{\mathsf{CSmProj}(k)}(S,S')$. We define:
\[
\begin{tabular}{cccc}
$(\omega_{M}(c))_r$  : & $\mathrm{Hom}_{\mathcal{M}(k)}(h(S),M(r))$ &  $\longrightarrow$  & $\mathrm{Hom}_{\mathcal{M}(k)}(h(S),M(r +r'))$ \\
 & $a$ & $\longmapsto$ & $(a \otimes \mathbbm{1}(r')) \circ c$.
\end{tabular}
\]
Here $c$ in the composition $(a \otimes \mathbbm{1}(r')) \circ c$ is viewed as a morphism $h(S') \longrightarrow h(S)(r')$ in $\mathcal{M}(k)$, and $\mathbbm{1}(r')$ denotes the $r'$-th Tate twist.
\end{definition}

\begin{remark} 

     Note that the sum of two correspondences of different degrees does not send a homogeneous element on a homogeneous element.
\end{remark}

\begin{definition}\label{def_homogeneous_map}
Let $M$ and $N$ be pure Chow motives, and let $\omega_M$, $\omega_N$ be the functors
defined in Definition~\ref{def_omega}.
\begin{enumerate}
    \item[\textup{(i)}] We define $\Hom(\omega_M,\omega_N)$ to be the set of natural
    transformations from $\omega_M$ to $\omega_N$.
    \item[\textup{(ii)}] A natural transformation $F \colon \omega_M \to \omega_N$ is
    said to be \emph{homogeneous of degree $d$} if, for every smooth projective
    $k$-scheme $S$, the induced map of graded modules
    \[
        F(S) \colon \omega_M(S) \longrightarrow \omega_N(S)
    \]
    is homogeneous of degree $d$, that is, for every $r \in \mathbb{Z}$, it sends
    $\omega_M(S)_r$ to $\omega_N(S)_{r+d}$.
    \item[\textup{(iii)}] We define $\Hom_{\mathrm{deg}\,0}(\omega_M,\omega_N)$ to be
    the subgroup of $\Hom(\omega_M,\omega_N)$ consisting of homogeneous natural
    transformations of degree $0$.
\end{enumerate}
\end{definition}

Given a morphism of motives $f \colon M \to N$, we associate a natural transformation of degree $0$,
\[
\Omega(f) \colon \omega_M \longrightarrow \omega_N,
\]
defined as follows: for every smooth projective $k$-scheme $S$ and every $r \in \mathbb{Z}$,
\begin{equation}\label{def_omega_function}
\begin{array}{cccc}
(\Omega(f)(S))_r \; : & \Hom_{\mathcal{M}(k)}(h(S),M(r)) & \longrightarrow & \Hom_{\mathcal{M}(k)}(h(S),N(r)) \\
& c & \longmapsto & (f \otimes \mathbbm{1}(r)) \circ c.
\end{array}
\end{equation}

The Yoneda lemma does not apply directly to the functor $\omega$, since the latter is defined on the category $\mathsf{CSmProj}(k)^{\mathrm{op}}$ rather than on $\mathcal{M}(k)$. The following result can be viewed as a variant of the Yoneda lemma.

\begin{proposition}[Manin Identity Principle] \label{Manin_Principle} 
Let $k$ be a field, and let $M$ and $N$ be pure Chow motives over $k$. Then the map $\Omega$ defined in~\eqref{def_omega_function}
\[
\Omega : \Hom_{\mathcal{M}(k)}(M,N) \longrightarrow \Hom_{\deg 0}(\omega_M, \omega_N)
\]
is an isomorphism, where $\Hom_{\deg 0}(\omega_M, \omega_N)$ is defined in Definition~\ref{def_homogeneous_map}.
\end{proposition}
\begin{proof} 
For any $r \in \mathbb{Z}$, tensoring with the motive $\mathbbm{1}(r)$ gives a bijection between maps $M \longrightarrow N$ and maps $M(r) \longrightarrow N(r)$. Moreover, giving a map of degree zero $\omega_{M} \longrightarrow \omega_{N}$ is the same as giving a map of degree zero $\omega_{M(r)} \longrightarrow \omega_{N(r)}$. Since these two identifications commute with the maps $\Omega$, we may assume that $M = (X, p, 0)$ 
for some projective scheme $X$ over $k$ and some projector $p \in \Corr^0(X, X)$. Then $\Hom_{\mathcal{M}(S)}(M, N)$ 
is the image of the endomorphism
$F \longmapsto F \circ p$
on $\Hom_{\mathcal{M}(S)}(h(X), N)$, and similarly 
$\Hom(\omega_M, \omega_N)$ 
is the image of the endomorphism
$f \longmapsto f \circ \Omega(p)$
on $\Hom(\omega_{h(X)}, \omega_N)$. Thus we are reduced to the case $M = h(X)$, where $X$ is a smooth projective scheme over $k$. We may further assume that $X$ has pure dimension $d$. Given a natural transformation $F : \omega_M \longrightarrow \omega_N$ of degree zero, we can associate to $F$ a morphism  
$f := F(\id_M)$, where we view $\id_M$ as an element of $(\omega_M(M))_0$. Then $f$ lies in $(\omega_N(M))_0 = \Hom_{\mathcal{M}(S)}(M, N)$ and since $F$ is functorial with respect to correspondences by definition, the image of $\id_M$ entirely determines $F$. Indeed, let $S \in \mathsf{SmProj}(k)$, let $r \in \mathbb{Z}$, and let 
$x \in (\omega_M(S))_r = \mathrm{Hom}_{\mathcal{M}(k)}(h(S), M(r))=\mathrm{Corr}^r(S,X)$. Then $x$ is a correspondence of degree $r$ from $S$ to $M$, hence may be viewed as a morphism in $\mathsf{CSmProj}(k)$. By Definition~\ref{def_omega}, the map 
$\omega_M(x) : \omega_M(M) \to \omega_M(S)$
sends $\mathrm{id}_M$ to $x$. It follows that
\[F(x) = F(\omega_M(x)(\mathrm{id}_M)) 
= \omega_N(x)(F(\mathrm{id}_M)) 
= (f \otimes \mathbbm{1}(r)) \circ x 
= \Omega(f)(x).\]
Here the third equality follows from the definition of the functoriality of $\omega_N$ (see ~\ref{def_omega}), and the last equality from the definition of $\Omega(f)$ in \ref{def_omega_function}. Hence $\Omega(f) = F$. This proves surjectivity of $\omega$. Moreover, as $\omega(f)(M)(\id_M)=f$, the map $\omega$ is injective.
\end{proof}

The following is a consequence of Proposition \ref{Manin_Principle}.

\begin{corollary}\label{exact_Manin_Principle}
Let $k$ be a field and let $M, M', M'' \in \mathcal{M}(k)$. Let
$f \colon M \longrightarrow M''$ and $i \colon M' \longrightarrow M$ be morphisms
of pure Chow motives.
\begin{enumerate}
    \item[\textup{(i)}] The morphism $i \colon M' \to M$ is a monomorphism in
    $\mathcal{M}(k)$ if and only if the induced morphism
    \[
    \omega_{M'} \longrightarrow \omega_M
    \]
    is an injective natural transformation of functors from
    $\mathsf{CSmProj}(k)^{\mathrm{op}}$ to $\mathsf{GrAb}$.

    \item[\textup{(ii)}] The morphism $i \colon M' \to M$ is the kernel of
    $f \colon M \to M''$ in $\mathcal{M}(k)$ if and only if the induced sequence
    \[
    0 \longrightarrow \omega_{M'} \longrightarrow \omega_M \longrightarrow \omega_{M''}
    \]
    is an exact sequence of functors from $\mathsf{CSmProj}(k)^{\mathrm{op}}$ to
    $\mathsf{GrAb}$.
\end{enumerate}
\end{corollary}

\begin{proof}
\textup{(i)} By Proposition~\ref{Manin_Principle}, the functor $\omega$ is fully
faithful, so $i$ is a monomorphism in $\mathcal{M}(k)$ if and only if
$\omega(i) \colon \omega_{M'} \to \omega_M$ is a monomorphism in the of functor
from $\mathsf{CSmProj}(k)^{\mathrm{opp}}$ to $\mathsf{GrAb}$. Since monomorphisms in $\mathsf{GrAb}$ are exactly injective maps, this is
equivalent to $\omega(i)$ being an injective natural transformation.

\noindent \textup{(ii)} Suppose first that $i \colon M' \to M$ is the kernel of
$f \colon M \to M''$ in $\mathcal{M}(k)$. Then $M'(r)$ is the kernel of
$M(r) \to M''(r)$ for all $r \in \mathbb{Z}$, and hence $\omega_{M'}$ is the kernel
of $\omega_M \to \omega_{M''}$, giving exactness of the sequence.
Conversely, suppose that $\omega_{M'}$ is the kernel of
$\omega_M \to \omega_{M''}$. Let $N \in \mathcal{M}(k)$ and let
$g \colon N \to M$ be a morphism such that $f \circ g = 0$. Then the composite
$\Omega(f) \circ \Omega(g) = 0$, so $\Omega(g)$ factors as a natural transformation
$\omega_N \to \omega_{M'}$ of degree zero. By Proposition~\ref{Manin_Principle},
this implies that $g$ factors through a morphism $N \to M'$ in $\mathcal{M}(k)$,
which is what we wanted to show.
\end{proof}

\subsubsection*{A monomorphism from $h(J)^{\vee}$ to $h(C^{[\bullet]})^{\vee} $}

To show that the result of Moonen and Polishchuk lifts to motives, we first prove that the maps $\widetilde{s}_S$ are functorial in $S$.
Recall that ${\mathbb{K}}_S$ is defined as:
\[ {\mathbb{K}}_S:= \text{Ker} \Big( P_{0,1}([{(p_0)}_S])\Big) \cap \bigcap_{n \geqslant 1} \text{Ker}\Big(P_{0,1}(C\times_k S)^{[n]}\Big) ~ \subset ~ \CH_*(C^{[\bullet]}/S)\]
where the operators $P_{0,1}$'s are defined in \eqref{Operators}.

\begin{lemma}\label{factorization_K}
Let $k$ be a field, and let $S$ and $S'$ be smooth projective schemes over $k$, with $c \in \Corr(S,S')$. Let $C$ be as in Notation~\ref{Notation_1} and $\mathbb{K}_S$ as in Lemma~\ref{K_S}. 
Then the map
\[
\CH_*(C^{[\bullet]}/S) \longrightarrow \CH_*(C^{[\bullet]}/S'), \quad x \longmapsto c \circ x
\]
sends $\mathbb{K}_S$ into $\mathbb{K}_{S'}$.
\end{lemma}
\begin{proof}
Since $\mathbb{K}_S$ is a submodule generated by its homogeneous elements with respect to the $C^{[\bullet]}$-grading on $\CH_*(C^{[\bullet]}/S)$, it suffices to prove that if $\alpha \in \mathbb{K}_S \cap \CH_*(C^{[m]}\times_k S)$ where $m$ is a non-negative integer and $c \in \CH^*(S\times_k S')$ then $\alpha \circ c \in \mathbb{K}_{S'}$. 
 
 By Lemma \ref{intersection_theory_lemma}, we have:
 \begin{itemize}
     \item $(i_m)_{S'}^*(c \circ \alpha) = c \circ (i_m)_S^*(\alpha) = 0$;
     \item $(\text{pr}_2)_*(\alpha_{m',m-m'})_{S'}^*(c \circ \alpha) =  c \circ (\text{pr}_2)_*(\alpha_{m',m-m'})_S^*(\alpha) = 0$.
 \end{itemize}
    Thus, by the definition of the operators given in~\eqref{Operators}, we conclude that $\alpha \in \mathbb{K}_{S'} $.
\end{proof}

\begin{proposition} \label{functorialite_s} Let $k$ be a field. Let $S$ and $S'$ be smooth projective schemes, and let $c \in \Corr(S,S')$. Let $C$, $J$, $(\widetilde{\sigma}_{S})_*$, and $(\widetilde{\sigma}_{S'})_*$ be as in Notation \ref{Notation_1}, and let $\mathbb{K}_S$ be as in Lemma \ref{K_S}. We then have the following commutative diagram:
\begin{equation}
\begin{tikzcd}[sep=huge]
{\mathbb{K}}_S \arrow[r, hookrightarrow] \arrow[d,"c \circ "']& \mathrm{CH}_*(C^{[\bullet]}/S) \arrow[r,"(\widetilde{\sigma}_{S})_*"] \arrow[d, "c \circ "'] & \mathrm{CH}_*(J \times_k S) \arrow[d,"c \circ "] \\
{\mathbb{K}}_{S'} \arrow[r, hookrightarrow] & \mathrm{CH}_*(C^{[\bullet]}/S') \arrow[r,"(\widetilde{\sigma}_{S'})_*"] & \mathrm{CH}_*(J \times_k S')
\end{tikzcd}
\end{equation}
and the section $\widetilde{s}_S$ is functorial in $S$ with respect to correspondences. 
\end{proposition}

\begin{proof}  The commutativity of the first square follows from Lemma~\ref{factorization_K}.
The commutativity of the second square follows directly from Lemma \ref{intersection_theory_lemma}, together with the identity
\[(\widetilde{\sigma}_{S})_* = \bigoplus_{m \in \mathbb{N}} \left((\sigma_{m})_{S}\right)_*.\]
Since the section  $\widetilde{s}_S$ is defined as the composition of the inverse of 
\[ \mathbb{K}_S \longhookrightarrow \CH_*(C^{[\bullet]}/S) \xrightarrow[]{(\widetilde{\sigma}_S)_*} \CH_*(J \times_k S) \] 
with the inclusion 
\[\mathbb{K}_S \longhookrightarrow \CH_*(C^{[\bullet]}/S), \]
both functorial in $S$ with respect to correspondences as shown above, we conclude that $\widetilde{s}_S$ is functorial in $S$.
\end{proof}

We now give an explicit expression of the functor $\omega_{h(Y)^{\vee}}$, where $Y$ is a smooth projective scheme of dimension $d'$. 
For any smooth projective scheme $Y$ of dimension $d'$, we have:
\begin{align}
\begin{split}
    \omega_{h(Y)^{\vee}}(S)& = \displaystyle \bigoplus_{r \in \mathbb{Z}} \mathrm{Hom}_{\mathcal{M}(k)}(h(S),h(Y)^{\vee}(r))\\
    &=\bigoplus_{r \in \mathbb{Z}} \mathrm{Hom}_{\mathcal{M}(k)}(h(S),h(Y)(d'+r))\\
    &=\bigoplus_{r \in \mathbb{Z}} \CH^{\mathrm{dim}\,S+d'+r}(S\times_k Y)\\
    &=\bigoplus_{r \in \mathbb{Z}} \CH_{-r}(S\times_k Y)\\
    &=\bigoplus_{r \in \mathbb{Z}} \CH_{-r}(Y \times_k S)
\end{split}
\end{align}

\noindent The first equality is the definition of $\omega_{h(Y)^{\vee}}$ (see Definition ~\ref{def_omega}), the second follows from the equality $h(Y)^{\vee} = h(Y)(d')$ (see~\cite[\S4.1]{André2004}), the third one follows from the definition of morphisms in $\mathcal{M}(k)$ and the last one is the equality induced by the swap map $S \times_k Y \longrightarrow Y \times S$.

\begin{remark} \label{remark_omega_dual} Let $k$ be a field. Let $S$ and $S'$ be smooth projective schemes over $k$. If $c \in \Corr(S',S)$, then the map
\[
\omega_{h(Y)^{\vee}}(S)=\bigoplus_{r \in \mathbb{Z}} \CH_{-r}(Y\times_k S) \longrightarrow \omega_{h(Y)^{\vee}}(S')=\bigoplus_{r \in \mathbb{Z}} \CH_{-r}(Y\times_k S')
\]
described in Remark~\ref{def_omega} is given by
\[
a \longmapsto {}^t c \circ a ,
\]
where ${}^t c \in \Corr(S, S')$ denotes the transpose of $c$. This follows from \cite[Prop.~16.1.1(b)]{Fulton1998}.
\end{remark}

\begin{lemma}\label{borne_grading}
   Let $k$ be a field. Let $S$, $C$, and $J$ be as in Notation~\ref{Notation_1}. Then the grading on $\mathrm{CH}_*(J \times_k S)$ is concentrated in degrees less than or equal to $3 \dim J$, that is,
$\widetilde{s}_S\bigl(\mathrm{CH}_*(J \times_k S)\bigr) \subset \bigoplus_{m=0}^{3\dim J} \mathrm{CH}_*\bigl(C^{[m]} \times_k S\bigr)$.
\end{lemma}
\begin{proof}
Let $g = \dim J$. By the definition of the grading, it is enough to prove that
\[\widetilde{s}_S\bigl(\mathrm{CH}_*(J \times_k S)\bigr)
~\subset~
\bigoplus_{m=0}^{3g} \mathrm{CH}_*\bigl(C^{[m]} \times_k S\bigr).
\]
Let $c \in \mathrm{CH}_*(J \times_k S)$. By Proposition~\ref{functorialite_s}, if $\Delta_J$ denotes the diagonal in $\mathrm{CH}_g(J \times_k J)$, then $\widetilde{s}_S(c)
= \widetilde{s}_S(c \circ \Delta_J)
= c \circ \widetilde{s}_J(\Delta_J)$.
Moreover, $\widetilde{s}_J(\Delta_J) \in \bigoplus_{m=0}^{3g} \mathrm{CH}_*\bigl(C^{[m]} \times_k J\bigr)$ by definition of $\widetilde{s}_J$ (see \ref{s_S_def}) and
\[
c \circ \left( \bigoplus_{m=0}^{3g} \mathrm{CH}_*\bigl(C^{[m]} \times_k J\bigr) \right)
~\subset~
\bigoplus_{m=0}^{3g} \mathrm{CH}_*\bigl(C^{[m]} \times_k S\bigr).
\]
This implies that $\widetilde{s}_S(c) \in \bigoplus_{m=0}^{3g} \mathrm{CH}_*\bigl(C^{[m]} \times_k S\bigr)$
\end{proof}

\begin{proposition}\label{monomorphism_3g}
Let $k$ be a field and let $C$ and $J$ be as in Notation~\ref{Notation_1}.
There is a monomorphism
\[
h(J)^{\vee} \longhookrightarrow \bigoplus_{m=0}^{3\dim J} h(C^{[m]})^{\vee}
\]
in $\mathcal{M}(k)$.
\end{proposition}

\begin{proof}
Let $g = \dim J$. By Proposition~\ref{Manin_Principle} and
Corollary~\ref{exact_Manin_Principle}, constructing a monomorphism
$h(J)^{\vee} \hookrightarrow \bigoplus_{m=0}^{3g} h(C^{[m]})^{\vee}$
is equivalent to constructing a degree-zero injection
$\omega_{h(J)^{\vee}} \hookrightarrow \bigoplus_{m=0}^{3g} \omega_{h(C^{[m]})^{\vee}}$.
By Remark~\ref{remark_omega_dual}, this amounts to constructing, for every smooth
projective $k$-scheme $S$, a degree-zero injection
\[
\bigoplus_{r \in \mathbb{Z}} \mathrm{CH}_{-r}(J \times_k S)
\longhookrightarrow
\bigoplus_{r \in \mathbb{Z}} \bigoplus_{m=0}^{3g} \mathrm{CH}_{-r}(C^{[m]} \times_k S),
\]
with respect to the grading induced by $r$, functorial in $S$ with respect to
correspondences. The morphism $\widetilde{s}_S$ provides such a map. By Lemma~\ref{borne_grading}
it induces a map between the above direct sums, it is functorial with respect to
correspondences by Proposition~\ref{functorialite_s}, and injective of degree zero
by Proposition~\ref{s_S}.
\end{proof}

\begin{remark}
In~\cite[Thm.~3.4]{MoonenPolishchuk2010}, the authors prove a stronger result.  
Since we will not need the full strength of their theorem, 
we only stated a weaker version, adding the factorization through the component $\bigoplus_{m=0}^{3g} h(C^{[m]})^{\vee}$.
\end{remark}

\subsubsection*{A filtration on $h(J)^{\vee}$}

Now that we have defined a monomorphism $h(J)^{\vee} \longhookrightarrow \bigoplus_{m=0}^{3g} h(C^{[m]})^{\vee}$, we would like to define a filtration 
\[
\cdots  \subset ({h(J)}^{\vee})^{\geqslant p+1}  \subset ({h(J)}^{\vee})^{\geqslant p}  \subset \cdots  \subset ({h(J)}^{\vee})^{\geqslant 0}={h(J)}^{\vee},
\]
such that $({h(J)}^{\vee})^{\geqslant p}$ is the kernel of the composite
\begin{equation}\label{i_p} 
i_p := h(J)^{\vee} \longhookrightarrow \bigoplus_{m=0}^{3g} h(C^{[m]})^{\vee} \longrightarrow  \bigoplus_{m=0}^{p} h(C^{[m]})^{\vee}
\end{equation}
for every $p\in \mathbb{N}$.

\begin{proposition}\label{defining_filtration}
Let $k$ be a field and let $C$ and $J$ be as in Notation~\ref{Notation_1}.
Let $p \in \mathbb{N}$. There exists an unique motive ${(h(J)^{\vee})}^{\geqslant p}$ in
$\mathcal{M}(k)$ whose associated functor $\omega_{{(h(J)^{\vee})}^{\geqslant p}}$
(see Definition~\ref{def_omega}) is given by
\[
\omega_{{(h(J)^{\vee})}^{\geqslant p}}(S)
\;=\; \bigoplus_{m \geqslant p} \mathrm{CH}_*^{[m]}(J \times_k S),
\]
and such that $({h(J)}^{\vee})^{\geqslant p}$ coincides with the kernel of the map
$i_p$ defined in~\eqref{i_p}.
\end{proposition}

\begin{proof}
We construct ${(h(J)^{\vee})}^{\geqslant p}$ as the kernel of a projector of
${h(J)}^{\vee}$. By the Manin Identity Principle, defining a projector $\pi_p$ of
${h(J)}^{\vee}$ amounts to constructing, for every $S \in \mathsf{SmProj}(k)$, a
degree-zero projection
\[
p_S \colon \mathrm{CH}_{-*}(J \times_k S) \longrightarrow \mathrm{CH}_{-*}(J \times_k S)
\]
functorial in $S$ with respect to correspondences.

Since the $C^{[\bullet]}$-grading on $\mathrm{CH}_*(C^{[\bullet]}/S)$ of
Definition~\ref{grading_def} and the section $\widetilde{s}_S$ are both functorial
in $S$ with respect to correspondences (see Proposition~\ref{functorialite_s}), the
induced grading
\[
\mathrm{CH}_{-*}(J \times_k S) = \bigoplus_{m \in \mathbb{N}} \mathrm{CH}_{-*}^{[m]}(J \times_k S)
\]
from Theorem~\ref{Theorem_4} is functorial in $S$ with respect to correspondence: for every
$c \in \mathrm{CH}_*(S' \times_k S)$ and $x \in \mathrm{CH}_{-*}^{[m]}(J \times_k S)$,
one has ${}^tc \circ x \in \mathrm{CH}_{-*}^{[m]}(J \times_k S')$. Hence the projection
\[
\pi_p \colon \mathrm{CH}_{-*}(J \times_k S)
\longrightarrow \bigoplus_{0 \leqslant m < p} \mathrm{CH}_{-*}^{[m]}(J \times_k S)
\longhookrightarrow \mathrm{CH}_{-*}(J \times_k S)
\]
is functorial in $S$ with respect to correspondences and of degree zero, meaning
$\pi_p\bigl(\mathrm{CH}_{-r}(J \times_k S)\bigr) \subset \mathrm{CH}_{-r}(J \times_k S)$
for every $r \in \mathbb{Z}$ as a composite of degree zero map. By the Manin Identity Principle \ref{Manin_Principle}, it therefore defines
a projector $\pi_p$ of ${h(J)}^{\vee}$.
We denote by $({h(J)}^{\vee})^{\geqslant p}$ the kernel of $\pi_p$. The description
of $\omega_{({h(J)}^{\vee})^{\geqslant p}}$ and the identification of
$({h(J)}^{\vee})^{\geqslant p}$ with the kernel of $i_p$ then follow from the
definition of $\pi_p$ and Corollary~\ref{exact_Manin_Principle}.
\end{proof}
We can now prove the desired result:

\begin{corollary}\label{filtration_motives}
Let $k$ be a field and let $C$ and $J$ be as in Notation~\ref{Notation_1}.
The motives $({h(J)}^{\vee})^{\geqslant p}$ defined in
Proposition~\ref{defining_filtration} fit into a filtration
\[
0 \subset ({h(J)}^{\vee})^{\geqslant 3g} \subset ({h(J)}^{\vee})^{\geqslant 3g-1}
\subset \cdots \subset ({h(J)}^{\vee})^{\geqslant 0} = {h(J)}^{\vee}
\]
in $\mathcal{M}(k)$, with $({h(J)}^{\vee})^{\geqslant p} = 0$ for all $p > 3g$.
Moreover, ${h(J)}^{\vee} = h(k) \oplus ({h(J)}^{\vee})^{\geqslant 1}$, and the map
$h(k) \to h(J)^{\vee}$ induced by $0 \in J$ is $\mathrm{id} \oplus 0$.
\end{corollary}

\begin{proof}
Since $i_p$ factors through $i_{p+1}$ for every $p \in \mathbb{N}$ and
$({h(J)}^{\vee})^{\geqslant p}$ is the kernel of $i_p$ by
Proposition~\ref{defining_filtration}, we have
$({h(J)}^{\vee})^{\geqslant p+1} \subset ({h(J)}^{\vee})^{\geqslant p}$.
Since $h(J)^{\vee} \hookrightarrow \bigoplus_{n=0}^{3g} h(C^{[n]})^{\vee}$
by Proposition~\ref{monomorphism_3g}, the map $i_p$ is a monomorphism for
$p \geqslant 3g$, and hence $({h(J)}^{\vee})^{\geqslant p} = 0$ for all $p > 3g$.
This yields the desired filtration. By \cite[Cor.~7.6]{MoonenPolishchuk2010},
$\mathrm{CH}_*^{[0]}(J \times_k S) = \mathrm{CH}_*(C^{[0]} \times_k S) = \mathrm{CH}_*(S)$
for all $S \in \mathsf{SmProj}(k)$, so
\[
\omega_{{h(J)}^{\vee}}
= \omega_{h(k)} \oplus \omega_{{(h(J)^{\vee})}^{\geqslant 1}},
\]
and the Manin Identity Principle gives the decomposition
${h(J)}^{\vee} = h(k) \oplus ({h(J)}^{\vee})^{\geqslant 1}$.
To identify the map $h(k) \to h(J)^{\vee}$ induced by $0 \in J$, consider the
diagram
\[
\begin{tikzcd}[column sep=2.5em]
\mathrm{CH}_*(S)
  \arrow[r] \arrow[rr, bend left=18, "i"]
& \bigoplus\limits_{m=0}^{3g} \mathrm{CH}_*^{[m]}(J \times_k S)
  \arrow[r, hook] \arrow[rr, bend right=18, "\mathrm{id}"']
& \bigoplus\limits_{m=0}^{3g} \mathrm{CH}_*(C^{[m]} \times_k S)
  \arrow[r, "(\widetilde{\sigma}_S)_*"]
& \bigoplus\limits_{m=0}^{3g} \mathrm{CH}_*^{[m]}(J \times_k S)
\end{tikzcd}
\]
where $i$ is the natural inclusion
$\mathrm{CH}_*(S) = \mathrm{CH}_*(C^{[0]} \times_k S) \hookrightarrow
\bigoplus_{m=0}^{3g} \mathrm{CH}_*(C^{[m]} \times_k S)$.
The composite $\mathrm{CH}_*(S) \to \bigoplus_{m=0}^{3g} \mathrm{CH}_*^{[m]}(J \times_k S)$
is therefore given by $[0]_*$, which shows that the map
$h(k) = h(k)^{\vee} \to h(J)^{\vee}$ is $\mathrm{id} \oplus 0$.
\end{proof}

\subsubsection*{The action of $[n]_*$ on the filtration}

We now study the action of $[n]_*$ on the filtration 
$(\omega_{({h(J)}^{\vee})^{\geqslant p}})_{0 \leqslant p \leqslant 3g}$. 
As before, our results on the motives will be deduced from the corresponding results 
\cite[Thm. 4]{MoonenPolishchuk2010} by applying the Manin Identity Principle (see Proposition~\ref{Manin_Principle}).

\begin{proposition}Let $k$ be a field. Let $C$ and $J$ be as in Notation \ref{Notation_1}. Let $n\in \mathbb{Z}$. The action of $[n]_*$ on $h(J)^{\vee}$ preserves the Chow motives $\left((h(J)^{\vee})^{\geqslant p}\right)_{p\in \mathbb{N}}$ defined in Proposition \ref{defining_filtration}, that is, for each $p\in\mathbb{N}$, $(h(J)^{\vee})^{\geqslant p}$ is stable under $[n]_*$.
in the category $\mathcal{M}(k)$.
Moreover, for every $p\in \mathbb{N}$, the endomorphism $[n]_*$ satisfies the following factorization:
\[
\begin{tikzcd}[row sep=huge, column sep=tiny]
{(h(J)^{\vee})}^{\geqslant p} \arrow[rr, "{[n]}_*-n^{p}"] \arrow[dr] & & {(h(J)^{\vee})}^{\geqslant p}\\
 & {(h(J)^{\vee})}^{\geqslant p+1} \arrow[hookrightarrow] {ur}
 \end{tikzcd}
\]
in the category $\mathcal{M}(k)$.
\end{proposition}
\begin{proof}
By the Manin Identity Principle (Proposition~\ref{Manin_Principle}), it suffices to prove these factorizations at the level of the $\omega$-functors (see Definition~\ref{def_omega}).  
According to Proposition~\ref{defining_filtration}, for every $p \in \mathbb{N}$, the functor $\omega_{({h(J)}^{\vee})^{\geqslant p}}$ is given by
\[
\begin{tabular}{cccc}
$\omega_{({h(J)}^{\vee})^{\geqslant p}} :$ & $\mathsf{SmProj}(k)$ & $\longrightarrow$ & $\mathsf{GrAb}$ \\
& $S$ & $\longmapsto$ & $\displaystyle \bigoplus_{m \geqslant p} \CH_*^{[m]}(J \times_k S)$.
\end{tabular}
\]
It follows that $({h(J)}^{\vee})^{\geqslant p}$ is stable under the action of $[n]_*$, since the filtration of Theorem~\ref{Theorem_4} is preserved by $[n]_*$.
The second factorization follows from the fact that $[n]_*$ acts as $n^m$ on $\mathrm{gr}_{\mathrm{Fil}}^m$ (see (ii) of Theorem~\ref{Theorem_4}).
\end{proof}
From the second factorization, we deduce by induction on $p$ the following result:

 \begin{corollary} \label{f_p_maps}
Let $k$ be a field. Let $C$ and $J$ be as in Notation \ref{Notation_1}. Let $n\in \mathbb{Z}$. Let $\left((h(J)^{\vee})^{\geqslant p}\right)_{p\in \mathbb{N}}$ be the Chow motives defined in Proposition \ref{defining_filtration}. For every $p \in \mathbb{N}$, there exist maps $f_p$ and $f_p'$ of the following form:
\[
f_p = [n]_*^p + \sum_{\substack{0 \leqslant i' \leqslant p-1 \\ 0 \leqslant j' \leqslant \frac{p(p+1)}{2}}} \varepsilon_{i',j'} [n]_*^{i'} n^{j'} 
\quad \text{with } \varepsilon_{i',j'} \in \mathbb{Z},
\]
\[
f_p' = [n]_*^p + \sum_{\substack{0 \leqslant i' \leqslant p-1 \\ 1 \leqslant j' \leqslant \frac{(p+1)(p+2)}{2}}} \varepsilon'_{i',j'} [n]_*^{i'} n^{j'} 
\quad \text{with } \varepsilon'_{i',j'} \in \mathbb{Z},
\]
such that the map $f_p : h(J)^{\vee} \longrightarrow h(J)^{\vee}$ factors through a map $\widetilde{f}_p : h(J)^{\vee} \longrightarrow (h(J)^{\vee})^{\geqslant p}$:
\[
\begin{tikzcd}[row sep=huge, column sep=small]
h(J)^{\vee} \arrow[rr, "f_p"] \arrow[dr, "\widetilde{f}_p"'] & & h(J)^{\vee} \\
& (h(J)^{\vee})^{\geqslant p} \arrow[hookrightarrow]{ur}
\end{tikzcd}
\]
and $f_p' : (h(J)^{\vee})^{\geqslant 1} \longrightarrow (h(J)^{\vee})^{\geqslant 1}$ factors through a map $\widetilde{f}_p' : (h(J)^{\vee})^{\geqslant 1} \longrightarrow (h(J)^{\vee})^{\geqslant p+1}$:
\[
\begin{tikzcd}[row sep=huge, column sep=small]
(h(J)^{\vee})^{\geqslant 1} \arrow[rr, "f_p'"] \arrow[dr, "\widetilde{f}_p'"'] & & (h(J)^{\vee})^{\geqslant 1} \\
& (h(J)^{\vee})^{\geqslant p+1} \arrow[hookrightarrow]{ur}
\end{tikzcd}.
\]
\end{corollary}

 \bigskip

\begin{corollary} \label{action_n_M}
Let $k$ be a field. Let $C$ and $J$ be as in Notation \ref{Notation_1}. Let ${(h(J)^{\vee})}^{\geqslant 1}$ be the Chow motive defined in Proposition \ref{defining_filtration}. Let $n\in \mathbb{Z}$.
\begin{enumerate}
    \item[\textup{(i)}]
    There exists an endomorphism $f'$ of ${(h(J)^{\vee})}^{\geqslant 1}$ such that the restriction of $[n]_*^{3g}$ to ${(h(J)^{\vee})}^{\geqslant 1}$ is equal to $n f'$.
    
    \item[\textup{(ii)}]
    There exists a decomposition of $h(J)$ of the form
    $h(J)=h(k) \oplus h'$, where $h'$ is an effective Chow motive. Moreover, $h'$ is stable under $[n]^*$ and the restriction of $([n]^*)^{3g}$ to $h'$ is of the form $n f$, where $f$ is an endomorphism of $h'$. Furthermore, the morphism
\[
h(J)= h(k) \oplus h' \longrightarrow h(k)
\]
induced by $0 \in J$ is the projection onto the first factor.
\end{enumerate}
\end{corollary}
\begin{proof}
Since ${(h(J)^{\vee})}^{\geqslant 3g+1}=0$ by Proposition \ref{defining_filtration}, the endomorphism $f'_{3g}$ of Corollary \ref{f_p_maps} is zero. This proves~(i).

\noindent Dualising the decomposition of Corollary \ref{filtration_motives} and using~(i), we obtain the desired decomposition $h(J)=h(k) \oplus h'$
in $\mathcal{M}(k)$. Since $h'$ is a direct summand of $h(J)$, it is an effective Chow motive. Hence the above decomposition already holds in $\mathcal{M}^{\mathrm{eff}}(k)$.
\end{proof}

\begin{remark}
The decomposition $h(J) = h(k) \oplus h'$ from (ii) of Corollary~\ref{action_n_M} may differ from the one induced by the point $0 \in J(k)$.
\end{remark}
\subsubsection*{Consequence for Voevodsky motives}

 In \cite[Prop.~2.1.4]{Voevodsky2000}, the author constructs an additive contravariant functor from the category of effective Chow motives $\mathcal{M}^{\mathsf{eff}}(k)$ to $\mathsf{DM}_{\mathrm{gm}}^{\mathsf{eff}}(k)$ such that the following diagram commutes:
\[
\begin{tikzcd}[sep=huge]
\mathsf{SmProj}(k) \arrow[r] \arrow[d,"\mathrm{M_{Chow}}"']  & \mathsf{Sm}(k) \arrow[d,"\mathrm{M_{gm}}"] \\
\mathcal{M}^{\mathsf{eff}}(k) \arrow[r] & \mathsf{DM}_{\mathrm{gm}}^{\mathsf{eff}}(k)
\end{tikzcd}
\]

Thus, Corollary~\ref{action_n_M}\,(ii), concerning the Chow motive $h(J)$, carries over to the category $\mathsf{DM_{gm}^{eff}}(k)$, and thus in $\mathsf{DM}(k)$.
\begin{proposition} \label{action_n_DM}
Let $k$ be a field. Let $C$ and $J$ be as in Notation~\ref{Notation_1}. Let $n\in \mathbb{Z}$. Then there exists a decomposition of $M(J)$ of the form
\[
M(J)=M(k) \oplus M'
\]
such that $M'$ is stable under the action of $\mathsf{M_{gm}}([n])$, and the restriction of $\mathsf{M_{gm}}([n])^{3g}$ to $M'$ is of the form $n f$, where $f$ is an endomorphism of $M'$. Moreover, the map
\[
M(k)  \longrightarrow M(J)= M(k) \oplus M'
\]
induced by $0 \in J$ is given by $\mathrm{id} \oplus 0$.
\end{proposition}

\section{Unramified cohomology and unramified Witt group of an abelian variety}

 In this section, we work over a field $k$ of characteristic $0$ and we will denote by $A$ an abelian variety over $k$ of dimension $g$.
 
In the preceding section, we associated to a Jacobian variety a filtration of its motive. We now aim to derive some consequences for the unramified cohomology and the unramified Witt group of an abelian variety. Although our filtration is only for Jacobian varieties, we will use Matsusaka's theorem --- which states that every abelian variety is a quotient of a Jacobian variety (see \cite[III.~§10]{milneAV}) --- together with Poincaré's complete reducibility theorem (\cite[Thm.~12.2]{EdixhovenVanDerGeerMoonenAV}) to extend our result to all abelian varieties, not just Jacobians of curves.

To this end, we will mainly rely on the following result:

\begin{proposition} \label{Matsusaka}
    Let $A$ be an abelian variety. Then there exists an étale isogeny of $A$ which factors through a Jacobian variety.
\end{proposition}

\begin{proof} By a theorem of Matsusaka (see \cite[III.~§10]{milneAV}), there exists a Jacobian variety such that $A$ is a quotient of $J$. Then, by Poincaré's complete reducibility theorem, (\cite[Thm.~12.2]{EdixhovenVanDerGeerMoonenAV}) There exists an abelian variety $B$ such that $A \times_k B$ is isogenous to $J$. In particular, there exists an integer $d \in \mathbb{N}$ such that the multiplication-by-$d$ map on $A \times_k B$ factors through $J$. It follows that the multiplication-by-$d$ map on $A$ factors through $J$. Since $\operatorname{char} k = 0$, this map is an étale isogeny.
\end{proof}

\subsection{Étale cohomology with finite coefficients}\label{section_etale_coho_finite_coeff}

Before dealing with the case of unramified cohomology and unramified Witt groups of abelian varieties, we begin by reviewing a similar result concerning the étale cohomology of abelian varieties with finite coefficients. The proof of this result is considerably more elementary and does not require the results from the work of Moonen and Polishchuk~\cite{MoonenPolishchuk2010} on the Chow groups of Jacobian varieties.

We first recall the following fact about the étale cohomology of abelian varieties.

\begin{proposition}\label{etale_finite_cohomology} 
    Let $k$ be an algebraically closed field of characteristic zero. Let $A$ be an abelian variety over $k$, and let $n \in \mathbb{N}$.
Then, for every $i \geqslant 1$, the cup product induces a natural isomorphism:
\[
\Lambda^i H^1(A, \mathbb{Z}/n) \cong H^i(A, \mathbb{Z}/n).
\]
\end{proposition}
\begin{proof}
\textit{First Step:}  
Let 
\[
0 \longrightarrow \mathbb{Z}/n' \longrightarrow \mathbb{Z}/n \longrightarrow \mathbb{Z}/n'' \longrightarrow 0
\]
be a short exact sequence, and suppose that the proposition holds for both $\mathbb{Z}/n'$ and $\mathbb{Z}/n''$. We aim to show that it also holds for $\mathbb{Z}/n$.

We always have a map
\[
\Lambda^i H^1(A, \mathbb{Z}/n) \longrightarrow H^i(A, \mathbb{Z}/n)
\]
given by the cup product. Therefore, it suffices to prove that this map is an isomorphism for all $i \geqslant 1$.

We proceed by induction on $i$ and show that the sequence
\begin{equation}\label{exact_sequence}
0 \longrightarrow H^i(A, \mathbb{Z}/n') \longrightarrow H^i(A, \mathbb{Z}/n) \longrightarrow H^i(A, \mathbb{Z}/n'') \longrightarrow 0
\end{equation}
is exact for every $i \geqslant 1$.

For $i=1$, note that 
\begin{equation}\label{dim_H^1}
H^1(A, \mathbb{Z}/n) = \operatorname{Hom}_{\text{cont}}(\pi_1(A), \mathbb{Z}/n) = \operatorname{Hom}_{\text{cont}}(\widehat{\mathbb{Z}}^{2g}, \mathbb{Z}/n) = (\mathbb{Z}/n)^{2g},
\end{equation}
for every $n \in \mathbb{N}^*$. Hence, the case $i=1$ follows.

Suppose that the sequence
\[
0 \longrightarrow H^i(A, \mathbb{Z}/n') \longrightarrow H^i(A, \mathbb{Z}/n) \longrightarrow H^i(A, \mathbb{Z}/n'') \longrightarrow 0
\]
is exact. By the long exact sequence in étale cohomology, we obtain
\[
0 \longrightarrow H^{i+1}(A, \mathbb{Z}/n') \longrightarrow H^{i+1}(A, \mathbb{Z}/n) \longrightarrow H^{i+1}(A, \mathbb{Z}/n'').
\]
The surjectivity for $i+1$ follows from the inductive hypothesis for $i=1$, the right-exactness of $\Lambda^{i+1}$, and the assumption that the proposition holds for $\mathbb{Z}/n''$.
Indeed, consider the following diagram:
\[
\begin{tikzcd}[row sep=huge, column sep=large]
\Lambda^{i+1} H^1(A, \mathbb{Z}/n) \arrow[r] \arrow[d] 
& \Lambda^{i+1} H^1(A, \mathbb{Z}/n') \arrow[d] \\
H^{i+1}(A, \mathbb{Z}/n) \arrow[r] 
& H^{i+1}(A, \mathbb{Z}/n'') 
\end{tikzcd}
\]
Since the top horizontal and right vertical maps are surjective, the bottom horizontal map is also surjective.
Moreover, we have the following commutative diagram of exact sequences:
\[ \begin{tikzcd}[row sep=huge, column sep=large] & \Lambda^i H^1(A, \mathbb{Z}/n') \arrow[r] \arrow[d] & \Lambda^i H^1(A, \mathbb{Z}/n) \arrow[r] \arrow[d] & \Lambda^i H^1(A, \mathbb{Z}/n'') \arrow[r] \arrow[d] & 0 \\ 0 \arrow[r] & H^i(A, \mathbb{Z}/n') \arrow[r] & H^i(A, \mathbb{Z}/n) \arrow[r] & H^i(A, \mathbb{Z}/n'') \end{tikzcd} \]
Since the first and last vertical arrows are isomorphisms by the inductive hypothesis, the middle vertical map is also an isomorphism by the Five Lemma.

\textit{Second Step:}  
By the first step, it suffices to prove the proposition in the case where $n = \ell$ is a prime number. In this situation, $H^*(A, \mathbb{Z}/\ell)$ carries the structure of a graded, associative, anticommutative bialgebra over the field $\mathbb{Z}/\ell$. The multiplication is given by the cup product, which is graded-commutative by \cite[1.18]{Milneetla}, while the comultiplication
\[
H^*(A, \mathbb{Z}/\ell) \longrightarrow H^*(A, \mathbb{Z}/\ell) \otimes_{\mathbb{Z}/\ell} H^*(A, \mathbb{Z}/\ell)
\]
arises from the group operation $m: A \times A \to A$ together with the Künneth formula for étale cohomology (see \cite[Prop.~7.4.11]{Fu2015}), and the counit is given by $[0]^*: H^*(A, \mathbb{Z}/\ell) \to \mathbb{Z}/\ell$.
By \eqref{dim_H^1}, we have $\dim_{\mathbb{Z}/\ell} H^1(A, \mathbb{Z}/\ell) = 2 \dim A$, and $H^i(A, \mathbb{Z}/\ell) = 0$ for all $i > 2 \dim A$ by \cite[Chap.~IV, Prop.~1.1]{Milne1986}. The result then follows from the Borel--Hopf structure theorem (see \cite[Thm.~15.2]{milneAV}), which describes the structure of graded connected anticommutative bialgebras.

\end{proof}

We deduce the following corollary:

\begin{corollary} \label{cancelling_etale_coho_A}Let $k$ be a field of characteristic $0$. Let $A$ be an abelian variety over $k$.
Let $\mu$ be a finite multiplicative group defines over $k$.  
Let $n \in \mathbb{N}^*$ be such that $n\mu = 0$. For every $i \geqslant 1$, there exists an integer $j \geqslant 1$ such that, for every $\alpha \in H^i(A,\mu)$, the class $(([n]^*)^j)(\alpha)$ is constant, that is, it comes from $H^i(k,\mu)$.

Moreover, if $k$ is algebraically closed, one can take $j=1$.
\end{corollary}

\begin{proof} 
First, we consider the case where $k$ is algebraically closed. In this situation, $\mu$ is a finite product of group schemes of the form $\mu_n$, i.e.
\[
\mu \cong \prod_i \mu_{n_i}
\quad \text{for some integers } n_i \in \mathbb{N}.
\]
We can suppose that $\mu=\mu_n$.
Since 
\[
H^1(A,\mathbb{Z}/n) \cong \operatorname{Pic}^0(A)[n],
\] 
the claim holds for $i=1$. By Proposition~\ref{etale_finite_cohomology}, the result then follows for every $i \geqslant 1$.

Now we turn to the general case. Let $k_s$ denote a separable closure of $k$.
We have the Hochschild–Serre spectral sequence
\[
E_2^{p,q} = H^p\!\bigl(k, H^q(A_{k_s},\mu)\bigr) \;\Longrightarrow\; H^{p+q}(A, \mu).
\]
Let $\alpha \in H^i(A, \mu)$.  
By the algebraically closed case, every term on the $E_2$-page can be killed after pulling back by $[n]^*$, except for the terms coming from $H^i\!\bigl(k, H^0(A_{k_s},\mu)\bigr)$.  
Since $A$ is proper, we have
\[
H^n\!\bigl(k, H^0(A_{k_s},\mu)\bigr) \cong H^n(k,\mu).
\]
Thus, after pulling back by some isogeny which is a power of $[n]$, we may assume that $\alpha$ comes from a constant element $\alpha' \in H^i(k,\mu)$.  
Moreover, this isogeny depends only on $i$ and not on $\alpha$.  
This gives the desired result.

        \end{proof}

\begin{corollary} \label{cancelling_etale_coho} Let $k$ be a field of characteristic $0$. Let $A$ be an abelian variety over $k$. Let $\mu$ be a finite multiplicative group defined over $k$. 
Let $U \subset A$ be an open subset that contains all points of codimension~$1$. 
Then, for $i = 1, 2$, there exists an integer $n$, an integer $j \geq 1$, 
and a morphism $f$ such that for every class 
$\alpha \in H^i(U, \mu)$, the element $(([n]^*)^j)(\alpha)$ is constant in 
$H^i((([n]^*)^j)^{-1}(U), \mu)$; that is, it comes from an element of $H^i(k, \mu)$.
\end{corollary}

\begin{proof}
Let $F := A \setminus U$ be the closed subscheme equipped with the reduced structure.  
Since $\operatorname{codim}(F) \geqslant 2$, and since $\mu$ satisfies the purity theorem (see \cite[Prop.~8.5.6]{Fu2015}), by the same argument as in the case of $\mu_m$ in \cite[Cor.~3.4.2]{ColliotThélène1995}, we obtain an isomorphism
\[
H^i(A, \mu) \longrightarrow H^i(U, \mu)
\]
for $i = 1, 2$.  
The result then follows directly from~\ref{cancelling_etale_coho_A}.
\end{proof}

We now state a more general result in the case of torsors under twisted constant group schemes over~$k$.

\begin{proposition} \label{cancelling_etale_coho_non_commutative} Let $k$ be a field of characteristic $0$
Let $G$ be a twisted constant group scheme over $k$ (that is not necessarily finite), and let $U \subset A$ be an open subset that contains all points of codimension $1$.  
Then there exists an étale isogeny $f$ of $A$ such that, for every $\alpha \in H^1(U,G)$, the class $f^*\alpha$ is constant in $H^1(f^{-1}(U),G)$; that is, it comes from an element of $H^1(k,G)$. 
In particular if $\alpha \in H^1(U,G)$ is a class such that $[0]^*\alpha \in H^1(k,G)$ is trivial, then, there exists an étale isogeny $f$ of $A$ such that the class $f^*\alpha$ is trivial in $H^1(f^{-1}(U),G)$.
\end{proposition}

\begin{proof}
By an extended version of the Zariski–Nagata purity theorem~\cite[Lem.~2.13~(3)]{Gille24}, 
we have an isomorphism
\[
H^{1}(A, G) \xrightarrow{\;\sim\;} H^{1}(U, G),
\]
so it suffices to prove the statement for $U = A$.

Let $\alpha \in H^{1}(A, G)$. 
By~\cite[Lem.~2.13~(1)]{Gille24}, the class \( \alpha \) is isotrivial. Hence, there exists a finite Galois extension \( k'/k \) and an integer \( n \in \mathbb{N}^{*} \) such that the composition
\[
A_{k'} \xlongrightarrow{p} A \xlongrightarrow{[n]} A
\]
trivializes $\alpha$. After pulling back by $[n]^*$, we may therefore assume that $\alpha$ is trivialized by the canonical projection $p : A_{k'} \longrightarrow A$.  

According to \cite[§2.2]{Gille2015}, we have the following exact sequence:
\[
1 \longrightarrow H^1(\mathrm{Gal}(k'/k), G(A_{k'})) 
\longrightarrow H^1(A,G) 
\longrightarrow H^1(A_{k'},G).
\]

Since $A$ is proper and $G$ is a twisted constant group scheme over $k$, we have
\[
H^1(\mathrm{Gal}(k'/k), G(A_{k'})) \;=\; H^1(\mathrm{Gal}(k'/k), G(k')).
\]

This shows that $\alpha$ comes from an element of $H^1(k,G)$ after pulling back by a suitable étale isogeny, as desired.
\end{proof}

From this proposition, we deduce the following corollary:

\begin{corollary}\label{torus_constant}Let $k$ be a field of characteristic $0$. Let $A$ be an abelian variety over $k$. Let $T$ be a torus defined over an open subset $U$ of $A$ that contains all points of codimension $1$. Then, there exists  an étale isogeny $f$ of $A$ such that $f^*T$ is a torus defined over $k$.   
\end{corollary}
\begin{proof}
    A torus of rank $n$ defined over $U$ is a form of $(\mathbb{G}_m)^n$ over $U$, and thus gives rise to a class
\[
    \alpha \in H^1(U, \mathrm{GL}_n(\mathbb{Z})).
\] 
    By Proposition~\ref{cancelling_etale_coho_non_commutative}, we can pull back $\alpha$ along an étale isogeny of $A$ so that it becomes a constant class. 
    This yields the desired result.
\end{proof}

\subsection{Unramified cohomology} \label{unramified_cohomology}

We would like to show that, for every $i \geqslant 1$, an element of $H^i_{\text{nr}}(A,\mathbb{Z}/n)$ becomes constant after pulling back via an étale isogeny of $A$. To prove this, we will use the Beilinson--Lichtenbaum conjecture (see \cite{SuslinVoevodsky2000}) which describes how unramified cohomology is related to motivic cohomology. 

\subsubsection*{Some general results on unramified cohomology}

\begin{notation}\label{Notation_H}
    Let $k$ be a field. Let $\mathcal{F}$ be an abelian sheaf on the étale site of $k$. We denote by $\mathcal{H}^j(\mathcal{F})$ the sheaf on the big Zariski site over $k$ associated to the presheaf:
    \[
    \begin{tabular}{cccc}
        $\mathsf{Sch}_k$ & $\longrightarrow$ & $\mathsf{Ab}$ \\
        $S$ & $\longmapsto$ & $H^j(S,\mathcal{F})$
    \end{tabular}
    \]
where $\mathsf{Sch}_k$ denotes the category of schemes over $k$ and $\mathsf{Ab}$ the category of abelian groups
\end{notation}

Let $X$ be a smooth, irreducible scheme of finite type and separated over the field $k$ and $\mathcal{F}$ is an étale sheaf over $X$ of the form $\mathbb{Z}/n(i)$ with $n \in \mathbb{N}$ and $i \in \mathbb{Z}$.

\noindent In \cite[Thm.~4.2]{BlochOgus1974}, Bloch and Ogus construct, for every $i \in \mathbb{N}$, the $i$-th \emph{Cousin complex} in étale cohomology
\begin{equation}0 \longrightarrow \bigoplus_{x \in X^{(0)}} H^i_x(X, \mathcal{F}) \longrightarrow
\bigoplus_{x \in X^{(1)}} H^{i+1}_x(X, \mathcal{F}) \longrightarrow \cdots \longrightarrow
\bigoplus_{x \in X^{(j)}} H^{i+j}_x(X, \mathcal{F}) \longrightarrow \cdots,
\end{equation}
with
$
H^i_x(X, \mathcal{F}) := \varinjlim_{x \in U} H^i_{\mathrm{\acute{e}t},\, \overline{\{x\}} \cap U}(U, \mathcal{F})$,
the inductive limit being taken over the open subsets  $U \subset X$ containing $x$.
Using Gabber's absolute purity theorem (see \cite[Exposé XVI, Thm. 3.1.1]{IllusieLaszloOrgogozo2014}), this complex can be rewritten as
\begin{equation}\label{Cousin_complexe}
\cdots \longrightarrow  H^{i}(k(X), \mathcal{F}) \longrightarrow
\bigoplus_{x \in X^{(1)}} H^{i-1}(k(x), \mathcal{F}(-1)) \longrightarrow \cdots \longrightarrow
\bigoplus_{x \in X^{(j)}} H^{i-j}(k(x), \mathcal{F}(-j)) \longrightarrow \cdots.    
\end{equation}
Bloch and Ogus show that the cohomology of this complex computes the Zariski cohomology of the sheaf $\mathcal{H}^p(\mathcal{F})$ on the small Zariski site of $X$, that is, of the sheaf associated to the presheaf
\[
U \longmapsto H^p(U, \mathcal{F}),
\]
where $U$ ranges over the small Zariski site of $X$.
\bigskip

\begin{definition}{\cite[p.11]{ColliotThélène1995}} Let $k$ be a field and $K/k$ a field extension. We define the \emph{unramified cohomology} of degree $i$ of $\mathbb{Z}/n$ as the group
\[H^i_{\mathrm{nr}}(K/k,\mathbb{Z}/n) := \bigcap_{R\in \mathcal{P}(K/k)} \operatorname{Im}(H^i(R,\mathbb{Z}/n) \longrightarrow H^i(K,\mathbb{Z}/n))\]
where $\mathcal{P}(K/k)$ is the set of discrete valuation rings containing $k$ and having fraction field $K$.
If $X$ is an integral $k$-scheme, we put
\[H^i_{\mathrm{nr}}(X, \mathbb{Z}/n):= H^i_{\mathrm{nr}}(K(X)/k,\mathbb{Z}/n).\]

\end{definition}

\begin{remark} \label{unramified_proper} If $X$ is proper, \cite[Thl. 4.1.1]{ColliotThélène1995} and the Bloch--Ogus theorem (see \cite[Thm. 4.2]{BlochOgus1974}) shows that
\[
H^i_{\mathrm{nr}}(X, \mathcal{F}) = H^0_{\mathrm{Zar}}(X, \mathcal{H}^i(\mathcal{F}))
\]
when $\mathcal{F}$ is of the form $\mathbb{Z}/n(i)$ with $n\in \mathbb{N}$ and $i \in \mathbb{Z}$.

\end{remark}
\subsubsection*{Beilinson--Lichtenbaum conjecture}

We now explain how unramified cohomology is related to motivic cohomology for Zariski sheaves of the form $\mathbb{Z}/{n}$.

In their article \cite{SuslinVoevodsky2000}, Suslin and Voevodsky showed that the Bloch--Kato conjecture is equivalent to the Beilinson--Lichtenbaum conjecture over a field of characteristic zero with coefficients in $\mathbb{Z}/n$. This implication was later established in positive characteristic by Geisser and Levine in \cite{Geisser2001BlochKato}. The Bloch--Kato conjecture, and hence the Beilinson--Lichtenbaum conjecture, was subsequently proved by Voevodsky, building on work of Rost (see \cite{Voevodsky2011}).

\begin{theorem}[Beilinson--Lichtenbaum conjecture]\label{Beilinson_Lichtenbaum}
Let $k$ be a field. Let $X$ be a smooth, irreducible, separated scheme of finite type over $k$. Let $n \in \mathbb{N}^*$. Then the natural map
\[
\mathbb{Z}/n(j) \xlongrightarrow{\sim} \tau_{\leqslant j} \mathbf{R}\alpha_*(\mathbb{Z}/n)
\]
is an isomorphism, where $(\alpha^*, \alpha_*)$ is the morphism of topoi induced by the natural morphism
from the big étale site of smooth schemes over $X$ to the big Zariski site of smooth schemes over $X$.
\end{theorem}



From the Beilinson–Lichtenbaum conjecture, we deduce the following result on motivic cohomology:

\begin{proposition}{\cite[Th.~1.3]{Totaro2003}}\label{coro_BL} Let $k$ be a field. Let $X$ be a smooth, irreducible, separated scheme of finite type over $k$. Let $n\in \mathbb{N}^*$ and $\mathcal{H}^j(\mathbb{Z}/n)$ be as defined in Notation \ref{Notation_H}. For each $j \geqslant 0$, there is a long exact sequence
\[
\cdots \to H^{i+j}_M(X, \mathbb{Z}/n(j-1)) 
\to H^{i+j}_M(X, \mathbb{Z}/n(j)) 
\to H^i_{\mathrm{Zar}}(X, \mathcal{H}^j(\mathbb{Z}/n)) 
\to H^{i+j+1}_M(X, \mathbb{Z}/n(j-1)) 
\to \cdots
\]
\end{proposition}

\begin{proof}
    There is a distinguished triangle in the derived category of Zariski sheaves on $X$:
    \[ \tau_{\leqslant j-1}\mathcal{R}\alpha_*(\mathbb{Z}/n) \longrightarrow \tau_{\leqslant j}\mathcal{R}\alpha_*(\mathbb{Z}/n) \longrightarrow \mathcal{H}^j(\mathbb{Z}/n)[-j].\]
    Hence, the Beilinson--Lichtenbaum conjecture (\ref{Beilinson_Lichtenbaum}) yields the desired result.
\end{proof}

\subsubsection*{Making classes in unramified cohomology constant by an étale isogeny}

 We now return to the particular case of abelian variety.

\begin{proposition} \label{cancelling_unramified_coho} 
Let $k$ be a field of characteristic zero. Let $A$ be an abelian variety over $k$. Let $n \in \mathbb{N}^*$. Then there exists an étale isogeny $f$ of $A$ such that, for every $\alpha \in H^i_{\mathrm{Zar}}\!\bigl(A, \mathcal{H}^j(\mathbb{Z}/n)\bigr)$, the pullback $f^*(\alpha)$ is constant, that is, it comes from an element of $H^i_{\mathrm{Zar}}\!\bigl(k, \mathcal{H}^j(\mathbb{Z}/n)\bigr)$. In particular, $f^*(\alpha)$ is trivial if $i \geqslant 1$.
\end{proposition}

 \begin{proof}
 \textit{First Step.} We begin by considering the case of a Jacobian variety $J$ of a smooth projective curve over $k$.  
Let $g$ denote its dimension.

We decompose
\[
M(J) = M(k) \oplus M',
\]
as in Proposition~\ref{action_n_DM}.
From Corollary~\ref{coro_BL}, we obtain the following exact sequence:
\[
H^{i+j}_M(J, \mathbb{Z}/n(j)) \longrightarrow 
H^i_{\mathrm{Zar}}(J, \mathcal{H}^j(\mathbb{Z}/n)) \longrightarrow 
H^{i+j+1}_M(J, \mathbb{Z}/n(j-1)).
\]

Let $\alpha \in H^i_{\mathrm{Zar}}(J, \mathcal{H}^j(\mathbb{Z}/n))$.  
By subtracting a constant element, we may assume that $0^*\alpha = 0$ in 
$H^i_{\mathrm{Zar}}(k, \mathcal{H}^j(\mathbb{Z}/n))$.  
Let $\alpha'$ denote its image in $H^{i+j+1}_M(J, \mathbb{Z}/n(j-1))$.

\noindent Since $\alpha'$ vanishes in $H^{i+j+1}_M(\Spec k, \mathbb{Z}/n(j-1))$, we have
\[
\alpha' \in \Hom_{\mathsf{DM}(k)}(M', \mathbb{Z}/n(j-1)[i+j+1]) 
    \subset H^{i+j+1}_M(J, \mathbb{Z}/n(j-1)).
\]
Thus, by Proposition~\ref{action_n_DM}, the pullback of $\alpha'$ via $\mathrm{M_{gm}}(n)^{3g}$ vanishes.
We deduce that, after pulling back $\alpha$ by a suitable power of $[n]^*$ (independent of $\alpha$), it lies in the image of $H^{i+j}_M(J, \mathbb{Z}/n(j))$. Let $\widetilde{\alpha}$ be a lift of $\alpha$ in $\Hom(M', \mathbb{Z}/n(j-1)[i+j])$.
Since
\[
H^{i+j}_M(J, \mathbb{Z}/n(j)) 
  = \Hom_{\mathsf{DM}(k)}(M', \mathbb{Z}/n(j-1)[i+j]) 
    \oplus \Hom_{_{\mathsf{DM}(k)}}(M(k), \mathbb{Z}/n(j-1)[i+j]),
\]
we may subtract a constant element from $\widetilde{\alpha}$ so that it lies in $\Hom_{\mathsf{DM}(k)}(M', \mathbb{Z}/n(j-1)[i+j])$. 
Thus, we can assume that $\widetilde{\alpha} \in \Hom(M', \mathbb{Z}/n(j-1)[i+j])$.
We can then iterate the preceding argument to split $\widetilde{\alpha}$.

Therefore, we conclude that, after pulling back $\alpha$ by a suitable power of $[n]^*$ (independent of $\alpha$), it becomes constant.
 
\bigskip

\textit{Second Step}. We now consider the case of an abelian variety.

By Proposition \ref{Matsusaka}, there exists a Jacobian variety $J$, together with morphisms
\[
f : A \longrightarrow J \quad \text{and} \quad f' : J \longrightarrow A
\]
such that $f' \circ f$ is an isogeny. Moreover, from the discussion above, there exists an integer $d \geqslant 1$ such that every element of $H^i_{\mathrm{Zar}}(J, \mathcal{H}^j(\mathbb{Z}/n))$ becomes constant after pulling it back by $([n]^*)^d$ .
Since the following diagram commutes:
\[
\begin{tikzcd}[sep=huge]
H^i_{\mathrm{Zar}}(A, \mathcal{H}^j(\mathbb{Z}/n)) \arrow[r,"{([ n ]^*)^d}"] \arrow[d,"(f')^*"'] &
H^i_{\mathrm{Zar}}(A, \mathcal{H}^j(\mathbb{Z}/n)) \arrow[d,"(f')^*"] \\
H^i_{\mathrm{Zar}}(J, \mathcal{H}^j(\mathbb{Z}/n)) \arrow[r,"{([ n ]^*)^d}"] \arrow[d,"f^*"'] &
H^i_{\mathrm{Zar}}(J, \mathcal{H}^j(\mathbb{Z}/n)) \arrow[d,"f^*"] \\
H^i_{\mathrm{Zar}}(A, \mathcal{H}^j(\mathbb{Z}/n)) \arrow[r,"{([ n ]^*)^d}"] &
H^i_{\mathrm{Zar}}(A, \mathcal{H}^j(\mathbb{Z}/n)).
\end{tikzcd}
\]
the map
\[
H^i_{\mathrm{Zar}}\!\bigl(A, \mathcal{H}^j(\mathbb{Z}/n)\bigr) 
\xrightarrow{(f' \circ f \circ [n]^d)^*} 
H^i_{\mathrm{Zar}}\!\bigl(A, \mathcal{H}^j(\mathbb{Z}/n)\bigr)
\]
sends every element of $H^i_{\mathrm{Zar}}(A, \mathcal{H}^j(\mathbb{Z}/n))$ to a constant class.
\end{proof} 

\begin{corollary} \label{cancelling_unramified_coho_closed}Let $k$ a field of characteristic zero. Let $A$ be an abelian variety over $k$. Let $n \in \mathbb{N}^* $. Then there exists an étale isogeny $f$ of $A$ such that the map :
\[
f^*: H^i_{\mathrm{Zar}}\!\bigl(A, \mathcal{H}^j(\mathbb{Z}/n)\bigr) \longrightarrow H^i_{\mathrm{Zar}}\!\bigl(A, \mathcal{H}^j(\mathbb{Z}/n)\bigr) 
\]
 is zero for every $i+j \geqslant 1$.   
\end{corollary}
\begin{proof}
By Proposition~\ref{canceling_unramified_witt_group}, there exists an
isogeny $f$ of $A$ such that for every integers $i,j$, for every class
$\alpha\in H^i_{\mathrm{Zar}}
\bigl(X,\mathcal{H}^j(\mathbb{Z}/n)\bigr)$,
the class $f^*\alpha$ is pulled back from a class on $\Spec k$.
Since $k$ is algebraically closed, we have
$H^i_{\mathrm{Zar}}
\bigl(\Spec k,\mathcal{H}^j(\mathbb{Z}/n)\bigr)=0$
whenever $i+j\ge1$. Therefore,
$f^*\alpha=0$ whenever $i+j\ge1$.
\end{proof}

\subsection{Unramified Witt group} \label{unramified_witt_group}

In this section, in the same way as for unramified cohomology, we show that for $i \geqslant 1$, the elements of the unramified Witt group $H^i(W_A)$ of an abelian variety $A$ become constant after being pulled back by a suitable étale isogeny. This result will, in fact, follow from the case of unramified cohomology, using a spectral sequence stated in \cite{Totaro2003}. From this, we will deduce further results on the derived Witt group of Balmer, using the Gersten spectral sequence (see \cite{BalmerWalter2002}).

\subsubsection*{The Gersten--Witt complex}

Let $X$ be a regular, irreducible scheme of finite type over~$k$.

The \emph{Gersten--Witt complex} $W_X$ of $X$, introduced by Balmer and Walter in \cite{BalmerWalter2002}, is a cochain complex of the form
\[
0 \longrightarrow W(k(X)) \longrightarrow \bigoplus_{\alpha \in X^{(1)}} W(k(x)) \longrightarrow \cdots \longrightarrow \bigoplus_{x \in X^{(n)}} W(k(x)) \longrightarrow \cdots,
\]
where $X^{(i)}$ denotes the set of codimension-$i$ points of $X$.

More precisely, this complex can be written as
\begin{equation}\label{Gersten_Witt_complex}
 0 \longrightarrow W(k(X)) \longrightarrow \bigoplus_{x \in X^{(1)}} W(k(x); \det(\mathfrak{m}_x/\mathfrak{m}_x^2)) \longrightarrow \cdots \longrightarrow \bigoplus_{x \in X^{(j)}} W(k(x); \ \det(\mathfrak{m}_x/\mathfrak{m}_x^2)) \longrightarrow \cdots,
\end{equation}
where $m_x$ is the maximal ideal of $\mathcal{O}_{X,x}$ for every $x\in X$, and $W(F;L)$ denotes the Witt group of quadratic forms over the field $F$ with values in a given line $L$, defined for example in \cite[Annexe E]{Fasel2008}.
The cohomology groups of this complex are denoted by $H^*(W_X)$.  The term $H^0(W_X)$ is also called the \emph{unramified Witt group} of $X$, and is denoted by $W_{\mathrm{nr}}(X)$.
The maps in this complex are called \emph{residue maps}. We recall their construction below in the non-oriented case.

\begin{proposition}[{\cite[Lem.~19.10]{Merkurjev2008}}]\label{residue_map}
Let $R$ be a discrete valuation ring (DVR) with fraction field $K$ and residue field $\kappa$, and let $\pi$ be a uniformizer of $R$.  
There exists a group homomorphism, depending on the choice of $\pi$,
\[
\partial_{\pi} \colon W(K) \longrightarrow W(\kappa),
\]
satisfying
\[
\partial_{\pi}\big( \langle u\pi^n \rangle \big) =
\begin{cases}
\langle \overline{u} \rangle & \text{if $n$ is odd}, \\[4pt]
0 & \text{if $n$ is even}.
\end{cases}
\]
\end{proposition}

These residue maps satisfy the following properties.

\begin{proposition}[{\cite[Lem.~19.14]{Merkurjev2008}}] \label{residue_maps_ideal}
Let $R$ be a discrete valuation ring (DVR) with fraction field $K$ and residue field $\kappa$, and let $\pi$ be a uniformizer of $R$. Let $I(K)$ and $I(\kappa)$ denotes the fundamental ideals of the Witt groups $W(K)$ and $W(\kappa)$, respectively.  
Then, for every integer $n \geqslant 1$, we have
\[
\partial_{\pi}\big(I^n(K)\big) \subset I^{n-1}(\kappa).
\]
\end{proposition}

\begin{lemma}{\cite[§2.1,§2.2]{Scharlau85}} \label{exacte_sequence_Witt_group} Let $R$ be a DVR with fraction field $K$ and residue field $\kappa$. If $2 \in R^{\times}$, then the following sequence  
\[ 0 \longrightarrow W(R) \longrightarrow W(K) \xrightarrow{\partial_{\pi}} W(\kappa) \longrightarrow 0\]
is an exact sequence of abelian groups.
\end{lemma}

We can also generalize the notion of unramified Witt groups to the powers of the fundamental ideal.

\begin{definition} For every irreducible scheme $X$ smooth over $k$, every DVR $R$ containing a field of characteristic different from $2$ with fraction field $K$, and every $i \in \mathbb{N}$, we define the abelian groups $I^i_{\mathrm{nr}}(R)$ and $I^i_{\mathrm{nr}}(X)$ as follows:
\begin{itemize}
    \item $I^i_{\mathrm{nr}}(R):=I^i(K) \cap W(R)$
    \item $I^i_{\mathrm{nr}}(X):=I^i(k(X)) \cap W_{\mathrm{nr}}(X)$.
\end{itemize}  
We will simply write $I_{\mathrm{nr}}(X)$ for $I^1_{\mathrm{nr}}(X)$.
\end{definition}

\subsubsection*{Link between unramified cohomology and unramified Witt groups}
In \cite[Thm.~1.1]{Totaro2003}, it is stated that there exists a spectral sequence, which can be obtained by filtering the Balmer--Walter complex.
Indeed, since the differential of the Balmer--Walter complex respects the filtration by powers of the fundamental ideal 
(see \cite[§9.2]{Fasel2008}), Massey’s method of exact couples yields the following spectral sequence, 
which converges whenever the cohomological $2$-dimension of the base field is finite.

\begin{proposition}{\cite[Thm. 1.1]{Totaro2003}}\label{Pardon_spectral_sequence}
   For any smooth, irreducible scheme $X$ over a perfect field $k$ of characteristic different from $2$ and of finite cohomological $2$-dimension, there exists a spectral sequence,
\[
E_2^{i,j} = H^{i+j}_{\mathrm{Zar}}(X, \mathcal{H}^j(\mathbb{F}_2))
\;\Longrightarrow\;
H^{i+j}(W_X),
\]
whose differentials have bidegree $(r, r-1)$ for $r \geqslant 2$. Moreover, the groups $H^i(X, \mathcal{H}^j)$ vanish unless $0 \leqslant i \leqslant j$.
\end{proposition}

The existence of this spectral sequence relies on the following theorem, which was originally conjectured by Milnor~\cite{Milnor70}. It was subsequently proved by Arason for $n = 3$~\cite{arason1975cohomologische}, by Rost for $n = 4$~\cite{Rost1986}, and in full generality by Orlov--Vishik--Voevodsky~\cite{Voevodsky_2003b, Orlov-Vishik-Voevodsky2007}.

\begin{theorem}{\cite[4.1]{Orlov-Vishik-Voevodsky2007}}\label{Milnor_invariants}
For every $i \in \mathbb{N}$, there exists a bijection
\[
e_i : I^i K / I^{i+1} K \xlongrightarrow{\sim} H^i(K, \mathbb{F}_2)
\]
which sends an $i$-fold Pfister form $\langle 1, -a_1 \rangle \otimes \cdots \otimes \langle 1, -a_i \rangle$ to the cup product $(a_1) \cdots (a_i)$.
\end{theorem}

We may combine this property with Proposition~\ref{cancelling_unramified_coho} to deduce the following corollary:

\begin{corollary}\label{canceling_unramified_witt_group} Let $k$ be a field of characteristic zero and of bounded cohomological dimension. Let $A$ be an abelian variety over $k$. For every $n \geqslant 0$, there exists an étale isogeny $f$ of $A$ such that the pullback map
\[
f^* : H^n(W_A) \longrightarrow H^n(W_A)
\]
sends every element of $H^n(W_A)$ to a constant class, that is, to an element coming from $H^n(W_k)$.
\end{corollary}
\begin{proof}
By Proposition~\ref{Pardon_spectral_sequence}, the group $H^n(W_A)$ admits a finite filtration whose graded pieces are quotients of subgroups of
\[
H^{i+j}_{\mathrm{Zar}}\!\left(A, \mathcal{H}^j(\mathbb{F}_2)\right), 
\qquad i+j = n, \; 0 \leqslant i \leqslant j.
\]
Since $\mathcal{H}^j(\mathbb{F}_2)=0$ for $j > 2\dim A + \operatorname{cd}(k)$, by \cite[Cor. VI 4.1]{Milneetla}, this filtration is finite. 
By iteratively applying Proposition~\ref{cancelling_unramified_coho}, and using the functoriality of the spectral sequence of Proposition~\ref{Pardon_spectral_sequence} with respect to pullbacks along smooth morphisms between smooth equidimensional $k$-schemes, we obtain the desired result.
\end{proof}

In the same way, from Corollary~\ref{cancelling_unramified_coho_closed}, we obtain the following corollary.

\begin{corollary} \label{canceling_unramified_witt_group_closed} Let $k$ be an algebraically closed field of characteristic zero.
Let $A$ be an abelian variety over the algebraically closed field $k$.
For every $n\geqslant 1$, there exists an étale isogeny $f$ of $A$ such that the pullback map
\[
f^* : H^n(W_A) \longrightarrow H^n(W_A)
\]
is zero.
\end{corollary}

\begin{remark} In the case of the Jacobian of a curve, the proof of Proposition~\ref{cancelling_unramified_coho} shows that $f$ can be taken to be a power of the multiplication-by-$2$ map.    
\end{remark}

In the case of an algebraically closed field, the unramified Witt group $W_{\mathrm{nr}}(X) = H^0(X, W_X)$ cannot be annihilated in the same way. Indeed, quadratic forms of \emph{odd rank} in $W_{\mathrm{nr}}(X)$ cannot be made to vanish by pulling them back along an étale isogeny. However, a similar result holds if we restrict to quadratic forms of \emph{even rank}. More generally, we have the following lemma. 

\begin{lemma}\label{injection_Witt_k(X)}
Let $k$ be a field of characteristic zero and let $X$ be a smooth integral $k$-scheme with a $0$-cycle of degree one. 
\begin{enumerate}
    \item[\textup{(i)}] Let $\mathcal{F}$ an abelian sheaf. Then the natural map $H^i(k,\mathcal{F}) \longrightarrow H^i(k(X),\mathcal{F})$ is injective.
    \item[\textup{(ii)}] Let $W(k)$ and $W(k(X))$ be the Witt groups of $k$ and $k(X)$, respectively, and let $I(k)$ and $I(k(X))$ denote the fundamental ideals of $W(k)$ and $W(k(X))$. Let $q$ be a quadratic form over $k$. Then its class $[q] \in W(k)$ lies in $I^n(k)$ if and only if the class of $q_{|k(X)}$ in $W(k(X))$ lies in $I^n(k(X))$.
\end{enumerate}
\end{lemma}

\begin{proof}
By the Chow moving lemma (see \cite[0B1U]{Stacks}), every non-empty open subscheme $U \subset X$ admits a zero-cycle of degree $1$. It follows that, for every open subscheme $U \subset X$, the natural map $H^i(k,\mathcal{F}) \to H^i(U,\mathcal{F})$ is injective. Indeed, let $U \subset X$ be an open subscheme and let $\sum_i n_i[x_i]$ be a zero-cycle of $U$ of degree $1$. The composite
\[
H^i(k,\mathcal{F}) \longrightarrow H^i(U,\mathcal{F}) \longrightarrow \bigoplus_i H^i(k(x_i),\mathcal{F}) \xlongrightarrow{\oplus_i n_i\,c_{k(x_i)/k}} H^i(k,\mathcal{F}),
\]
where $c_{k(x_i)/k}$ denotes the corestriction map in Galois cohomology defined in \cite[Appendix~2.2.2]{Serre2002}, is the identity on $H^i(k,\mathcal{F})$. This proves the injectivity. Passing to the limit over all such open subschemes $U$ (see \cite[Cor. 5.9.3]{Fu2015}), we obtain the desired result.

\noindent The second point follows from the Milnor invariant (see Theorem \ref{Milnor_invariants}) together with the first point. We only prove the converse, the other implication being immediate. Let $q$ be a quadratic form over $k$. If $[q_{|k(X)}] \in I^n(k(X))$, then an induction on $i$, using the Milnor invariant and the injectivity proved in (i), shows that $[q] \in I^i(k)$ for every $i \le n$. Hence $[q] \in I^n(k)$.
\end{proof}

\begin{proposition} \label{cancelling_unramified_Witt_group_zero}
Let $k$ be a field of bounded cohomological dimension $n$. Let $A$ be an abelian variety defined over $k$. Then there exists an étale isogeny $f$ of $A$ such that the map
\[
f^* : I^i_{\mathrm{nr}}(A) \longrightarrow I^i_{\mathrm{nr}}(A)
\]
is the zero map for every $i>n$.
\end{proposition}

\begin{proof}
From Corollary \ref{canceling_unramified_witt_group}, we know that there exists an étale isogeny
\[
f^* : I^i_{\mathrm{nr}}(A) \longrightarrow I^i_{\mathrm{nr}}(A)
\]
which sends every element of $I^i_{\mathrm{nr}}(A)$ to a constant class, that is, to a class coming from $W(k)$. By Lemma \ref{injection_Witt_k(X)}, such a constant class lies in $I^i(k)$. Since $k$ has bounded cohomological dimension $n$, we have $I^i(k)=0$ for $i>n$. Therefore, the claim follows.
\end{proof}

\subsubsection*{Witt groups of an abelian variety}

Let $X$ be a smooth irreducible scheme over $k$.

\noindent One can define the \emph{derived Witt groups} $W^n(X)$ in the sense of Balmer.  We give a brief overview of these groups following \cite{BalmerWalter2002}; for further details, we refer the reader to the original article.

The groups $W^n(X)$ depend only on $n \pmod{4}$. Moreover, the even derived Witt groups are identified as follows:
\[
W^{4n}(X) \cong W(X), \quad W^{4n+2}(X) \cong W^{-}(X),
\]
where $W(X)$ and $W^{-}(X)$ denote the classical Witt groups of symmetric and skew-symmetric vector bundles on $X$, respectively (see Balmer \cite[Thm.~4.7]{Balmer2001}). The odd derived Witt groups $W^{4n+1}(X)$ and $W^{4n+3}(X)$ are isomorphic to the Witt groups of formations, as described by Fernández-Carmena \cite{Fernandez-Carmena1987}, Pardon \cite{Pardon1984}, and Ranicki \cite{Ranicki1973}, \cite{Ranicki1989}.

We have the following theorem (see \cite[Thm. 7.2]{BalmerWalter2002}), which links the derived Witt groups $W^n(X)$ to the unramified Witt groups:

\begin{theorem} \label{Spectral_sequence_Witt} Let $k$ be a field. Let $X$ be a smooth irreducible scheme over $k$. There exists a spectral sequence, called the \emph{Gersten--Witt spectral sequence}, which abuts to $E^n = W^n(X)$ such that 
\[
E_1^{p,q} \neq 0 \text{ if and only if } 0 \leqslant p \leqslant \dim(X),
\]
 and such that the horizontal lines
\[
0 \longrightarrow E_1^{0,q} \xrightarrow{~d_1~} E_1^{1,q} \xrightarrow{~d_1~} E_1^{2,q} \xrightarrow{~d_1~} \cdots \xrightarrow{~d_1~} E_1^{\dim(X),q} \longrightarrow 0
\]
of the spectral sequence vanish for $q \not\equiv 0 \pmod{4}$, and coincide with the Gersten--Witt complex (\ref{Gersten_Witt_complex}) when $q \equiv 0 \pmod{4}$.
\end{theorem}

In the case of an abelian variety $A$, the preceding arguments yield the following result:

\begin{corollary} Let $k$ be a field of characteristic $0$ and of bounded cohomological dimension. Let $A$ be an abelian variety over $k$. For every $n \in \mathbb{Z}$, there exists an étale isogeny $f$ of $A$ such that the image of the map
\[
f^* : W^n(A) \longrightarrow W^n(A)
\]
consists of constant classes, that is, of elements coming from $W^n(k)$.
\end{corollary}

\begin{proof}
Theorem~\ref{Spectral_sequence_Witt} provides a spectral sequence
\[
E^{p,q}_2 \Longrightarrow W^n(A),
\]
where each term $E^{p,q}_2$ is either zero or $H^p(W_A)$ with $p + q \equiv p \pmod{4}$.
The result then follows from Proposition~\ref{canceling_unramified_witt_group} and from the fact that, by construction, the Gersten--Witt spectral sequence is functorial with respect to pullbacks along flat morphisms between regular irreducible schemes over $k$.

\end{proof}

Proceeding with the same argument, but using Proposition~\ref{canceling_unramified_witt_group_closed} instead, we obtain the following result:

\begin{corollary}
Let $k$ be an algebraically closed field of characteristic $0$.  
For every $n \not\equiv 0 \pmod{4}$, there exists an étale isogeny $f$ of $A$ such that the image of the map
\[
f^* : W^n(A) \longrightarrow W^n(A)
\]
is zero.
\end{corollary}

\section{A question of Colliot-Thélène and Iyer}

In their paper \cite{ColliotThélèneIyer2011}, the authors raise a question concerning smooth projective families of homogeneous spaces over an abelian variety.

\begin{question}{\cite[Question~3.4]{ColliotThélèneIyer2011}}\label{conjecture_originale}
    Let $A$ be an abelian variety over an algebraically closed field $k$. Let $X \longrightarrow A$ be a smooth projective family of homogeneous spaces of connected linear algebraic groups. Does there exist a finite étale map $B \longrightarrow A$ with $B$ coonected such that $X \times_A B \longrightarrow B$ admits a rational section?
\end{question}

By a \textit{family of homogeneous spaces of connected linear algebraic groups}, we mean that for every point $\Spec F \longrightarrow A$, where $F$ is a field, the fiber $X_F$ is a homogeneous space under a connected linear algebraic group defined over $F$. When the family is projective, any connected solvable group acts trivially on the fibers (see \cite[Thm.~17.1]{Milne_2017}). Thus, the family can be assumed, without loss of generality, to be a family of homogeneous spaces over connected semisimple groups.

\begin{remark} \label{Rem_strut_av}
Since $B \longrightarrow A$ is finite étale and $k$ is algebraically closed, we can endow $B$ with the structure of an abelian variety such that $B \longrightarrow A$ is a morphism of abelian varieties (specifically, an isogeny; see \cite[Chap.~IV, §18]{Mumford1970}).
\end{remark} 

In this section, we provide an affirmative answer to Question \ref{conjecture_originale} in characteristic zero for reductive groups whose root systems do not contain a factor of type $E_8$. 

We actually consider a slightly more general setting, allowing the smooth projective family of homogeneous spaces to be defined over an open subset $U \subset A$ containing all points of codimension $1$, and working not only over algebraically closed fields but also over fields of cohomological dimension at most $1$. To achieve this, we begin by addressing the question for torsors under reductive groups (avoiding $E_8$). We then treat the case of smooth projective families of homogeneous spaces which, thanks to a theorem of Demazure (see \cite[Prop.~4]{Demazure1977}), requires no further work.

In the following, $k$ will always denote a field, and $A$ an abelian variety defined over $k$.

\subsection{The question for torsors} \label{The_conjecture_for_torsors}

We now consider Question~\ref{conjecture_originale} in a slightly different setting, namely the case where $X \longrightarrow A$
is a $G$-torsor, with $G$ a connected reductive group scheme defined over the abelian variety $A$. Moreover, we no longer assume that the base field $k$ is algebraically closed. The Grothendieck--Serre conjecture, proved by Fedorov and Panin in~\cite{Federov_Panin_2015} for infinite fields and in~\cite{Panin2020} for finite fields, asserts that a generically trivial $G$-torsor over a regular local ring is already trivial in the Zariski topology. In other words, for such torsors, generic triviality is equivalent to Zariski-local triviality.
This naturally leads to the following question in our setting :

\begin{question} \label{conjecture_A}
Let $A$ be an abelian variety over a field $k$, and let $G$ be a linear algebraic group defined over $A$.  
Let $X \longrightarrow A$ be a $G$-torsor trivial at the origin.  
Does there exist an étale isogeny $f : A \longrightarrow A$ such that the pullback $f^*X \longrightarrow A$ is Zariski-locally trivial?
\end{question}

By Remark \ref{Rem_strut_av}, a positive answer to this question will also give a positive answer to Question \ref{conjecture_originale}.

A positive response was already given by Berkovich in~\cite[p.~182]{Berkovich1973} for the case of $\mathrm{PGL}_n$ over an algebraically closed field.

A slightly more general question is to consider torsors not defined over $A$ itself but over a sufficiently large open subset of $A$:

\begin{question} \label{conjecture}
Let $A$ be an abelian variety over a field $k$, and let $U \subset A$ be an open subset that contains all points of codimension $1$, and the origin  $0$ of $A$.  
Let $G$ be a linear algebraic group defined over $U$, and let $X \longrightarrow U$ be a $G$-torsor trivial at the origin.  
Does there exist an étale isogeny $f : A \longrightarrow A$ such that the pullback $f^*X \longrightarrow f^{-1}(U)$
is Zariski-locally trivial?
\end{question}

A first remark is that when the group $G$ is defined over $k$, a positive answer to Question~\ref{conjecture} implies a positive answer to the following question.

\begin{question}
Let $A$ be an abelian variety over a field $k$, and let $U \subset A$ be an open subset containing the origin $0$ and all points of codimension $1$. Let $G$ be a linear algebraic group defined over $k$, and let $X \longrightarrow A$ be a $G$-torsor. Does there exist an étale isogeny $f \colon A \longrightarrow A$ such that the pullback
\[
f^*X|_{f^{-1}(U)} \longrightarrow f^{-1}(U)
\]
is generically induced by a $G$-torsor over $\Spec k$? That is, Does there exist a dense open subscheme $V \subset f^{-1}(U)$ such that the restriction $(f^*X)|_V$ arises from a constant $G$-torsor?
\end{question}

The remainder of this subsection is devoted to providing an affirmative answer to Question~\ref{conjecture} for reductive groups whose root datum does not contain a factor of type $E_8$, in the case where $k$ is of characteristic zero.

In the following, unless otherwise specified, $k$ will be of characteristic zero.

\subsubsection*{The case of semisimple simply connected quasi-split groups}

We begin by considering the case of semisimple, simply connected, split groups defined over $k$. In the case where the abelian variety $A$ has dimension $2$ and the group $G$ is a $k$-semisimple group, the answer to Question~\ref{conjecture} is a special case of a Theorem in \cite{Jong2011}.

\begin{theorem}{\cite[Theorem~1.4]{Jong2011}}
\label{thm:serre_conj_II}
Let $K$ be an algebraically closed field of arbitrary characteristic, and let $K'/K$ be the function field of a surface over $K$.  
Let $G$ be a connected, semisimple, simply connected algebraic group defined over $K$.  
Then every $G$-torsor over $K'$ is trivial.
\end{theorem}

Since any $k$-semisimple group is a direct product of almost simple, simply connected groups (see \cite[(26.8)]{Knus_et_al_1998}), it suffices to treat the case of almost simple, simply connected groups. Such groups are classified by their Dynkin diagrams (\emph{loc. cit.}).

\medskip

\textit{$\bullet$ The cases $A_n$ and $C_n$}:

\medskip
We begin by giving a classification of quasi-split, semisimple, simply connected $k$-groups of type $A_n$ ($n \geqslant 1$) and $C_n$ ($n \geqslant 1$).

The split, semisimple, simply connected groups of type $A_n$ are isomorphic to $\mathrm{SL}_{n+1}$, and the split, semisimple, simply connected groups of type $C_n$ are isomorphic to $\mathrm{Sp}_{2n}$.  
Moreover, since groups of type $C_n$ admit no outer automorphisms, all quasi-split groups of this type are split.  
Semisimple, simply connected, quasi-split groups of type $A_n$ are classified by hermitian forms (see \cite[Table p.55]{Tits1966} for $A_n$ and see \cite[Table p.56]{Tits1966} for $C_n$).

Let $\mathrm{SU}(h,k')$ denotes the unique quasi-split $k$-form of $\mathrm{SL}_{n+1}$ associated with a quadratic étale $k$-algebra $k'$ and a hermitian form $h$ over $k'/k$ of rank $n+1$ entirely determined by $k'$. $h$ is of the form
\[
h =
\begin{cases}
z_1\overline{z_{2}} \perp \cdots \perp z_n\overline{z_{n+1}}, & \text{if } n+1 = 2m, \\
z_0\overline{z_0} \perp z_1\overline{z_{2}} \perp \cdots \perp z_{n-1}\overline{z_{n}}, & \text{if } n+1 = 2m+1,
\end{cases}
\]
for some orthogonal basis $(z_1,\ldots,z_{n+1})$ or $(z_0,z_1,\ldots,z_{n})$ if $n+1$ is even or odd, respectively.
According to \cite[Ch.~10, §1.1]{Scharlau85}, one can associate to $h$ a non-degenerate quadratic form
\[
q_h = h + \overline{h},
\]
called the \emph{trace form}.  

This gives rise to a diagram of algebraic groups
\begin{equation}\label{diagramSU}
\begin{tikzcd}[sep=huge]
 \mathrm{SU}(h,k') \arrow[d,hookrightarrow] \arrow[r,hookrightarrow,"i"] & \mathrm{SO}(q_h) \arrow[d,hookrightarrow]  \\
 \mathrm{U}(h,k') \arrow[r,hookrightarrow] & \mathrm{O}(q_h) 
\end{tikzcd}
\end{equation}

\begin{lemma}\label{injectivity} Let $k$ be a field of characteristic zero, and $k'/k$ a quadratic extension of $k$ and $h$ an hermitian form over $k'/k$. Let $i :  \mathrm{SU}(h,k') \longrightarrow \mathrm{SO}(q_h)$ describe above.
The map $i_*:H^1(K, \mathrm{SU}(h,k')) \longrightarrow H^1(K,\mathrm{SO}(q_h))$ is injective for every field extension $K/k$.
\end{lemma}
\begin{proof}It is enough to prove the case $K = k$. By a twisting argument, it is sufficient to show that the map has trivial kernel. Since the determinant maps are surjective, the maps induced in cohomology by the horizontal maps in Diagram~\ref{diagramSU} also have trivial kernels.
Moreover, a cocycle calculation shows that the map
\[
H^1(k, \mathrm{U}(h,k')) \longrightarrow H^1(k, \mathrm{O}(q_h))
\]
which sends a Hermitian form to its trace form coincides with the map induced by $i$. By \cite[Ch.~10, §1.1]{Scharlau85}, this map is injective. We deduce that
\[
i_* \colon H^1(k, \mathrm{SU}(h,k')) \longrightarrow H^1(k, \mathrm{SO}(q_h))
\]
is injective.
\end{proof}

\begin{proposition}\label{conjecture_A_C} Question \ref{conjecture} admits a positive answer for semisimple simply connected quasi-split groups of types $A_n$  ($n \geqslant 1$) and $C_n$  ($n \geqslant 3$) defined over a field $k$ of characteristic zero, provided that $k$ has finite cohomological dimension or contains a square root of $-1$.
\end{proposition}

\begin{proof}
If $k$ contains a square root of $-1$, we may assume that $A$, $U$, $G$, and $X$ are defined over a subfield of $k$ finitely generated over $\mathbb{Q}(\sqrt{-1})$. We may therefore reduce to the case where $k$ has finite cohomological dimension.

  Since $\mathrm{SL}_n$ and $\mathrm{Sp}_{2n}$ are special, that is, they have no non-trivial torsors over a field, the split case is trivial. Let $U \subset A$ be an open subset as in Question~\ref{conjecture}.  Let $\mathrm{U}(h,k')$ denotes the unitary group of type $^2A_{n-1}$ associated with a quadratic extension $k'/k$, where $h$ is the hermitian form of rank $n$ associated with this quadratic extension , and let $q_h$ denotes the trace form of $h$.
Let $\alpha \in H^1(U, \mathrm{SU}(h,k'))$ be the class corresponding to an $\mathrm{SU}(h,k')$-torsor over $U$, which is trivial at the origin.  
Denote by $h_{\alpha}$ the hermitian form associated with $\alpha$.
Under the composite
\[
H^1(U, \mathrm{SU}(h,k')) \longrightarrow H^1(U,SO(q_h)) \longrightarrow W(k(A)),
\]
the image of $\alpha$ corresponds to a quadratic form $q_{h_{\alpha}}$ such that the difference $q_{h_{\alpha}} - q_h$ lies in the first power of the fundamental ideal, that is,
\[
q_{h_{\alpha}} - q_h \in I(k(A)) \subset W(k(A)).
\]
Indeed, since $q_{h_{\alpha}}$ is a form of $q_h$, the two quadratic forms have the same rank. Moreover, since both $q_{h_{\alpha}}$ and $q_h$ are defined over $U$, which contains all points of codimension~$1$, their difference lies in the Witt group of the local ring at every codimension-one point of $A$; that is,
\[
q_{h_{\alpha}} - q_h \in W(\mathcal{O}_{A,x}) \quad \text{for all } x \in A^{(1)}.
\]
By Lemma~\ref{exacte_sequence_Witt_group}, it follows that
\[
q_{h_{\alpha}} - q_h \in W_{\mathrm{nr}}(k(A)) \cap I(k(A)) = I_{\mathrm{nr}}(k(A)).
\]

Then, by Proposition~\ref{canceling_unramified_witt_group}, there exists an étale isogeny $f \colon A \to A$ such that $f^*(q_{h_{\alpha}} - q_h)$ becomes constant.  
We also have the following specialization diagram at $0 \in U(k)$, given in \cite[Prop.~2.1.10]{ColliotThélène1995}:
\begin{equation} \label{Witt_diagram}
\begin{tikzcd}[sep=huge]
W(U) \arrow[r,"s_0"] \arrow[d] & W(k) \arrow[d] \\
W_{\mathrm{loc}}(U) \arrow[r,"\partial_{0}^*"] & W_{\mathrm{loc}}(k)
\end{tikzcd}
\end{equation}
where $W_{\mathrm{loc}}$ is the functor defined in \cite[p.11]{ColliotThélène1995}, the top horizontal map is the pullback map at $0 \in U(k)$ and the bottom horizontal map comes from the functoriality of the functor $W_{\mathrm{loc}}$.
The class of $f^*(q_{h_{\alpha}} - q_h)$ in $W_{\mathrm{nr}}(A)$ is constant; hence it is also constant in $W_{\mathrm{loc}}(A) \subset W_{\mathrm{nr}}(A)$.  
Since $\alpha \in H^1(U, \mathrm{SU}(h,k'))$ is trivial at the origin, (\ref{Witt_diagram}) implies the triviality of $f^*(q_{h_{\alpha}} - q_h)$, and hence the triviality of $f^*(i(\alpha))$ in $H^1(k(A), \mathrm{SO}(q_h))$.
The kernel of the natural map
\[
i: H^1(k(A), \mathrm{SU}(h,k')) \longrightarrow H^1(k(A), \mathrm{SO}(q_h))
\]
is trivial by lemma \ref{injectivity}.  
It follows that $f^*(\alpha)_{|{k(A)} = 0}$ in $H^1(k(A), \mathrm{SU}(h,k'))$.
This completes the proof.
\end{proof}

\textit{$\bullet$ The cases $B_n$ and $D_n$:}

\medskip

We begin by giving a classification of quasi-split, semisimple, simply connected $k$-groups of type $B_n$ ($n \geq 2$) and $D_n$ ($n \geq 4$).
The split, $k$-semisimple, simply connected groups of type $B_n$ and $D_n$ are of the form $\mathrm{Spin}(q_0)$,  
where $q_0$ is a hyperbolic quadratic form. Moreover, since groups of type $B_n$ have no outer automorphisms, all quasi-split groups of this type are split.  
For quasi-split groups of type $D_n$, they are isomorphic to $\mathrm{Spin}(q_0)$, where the quadratic form is given by
\[
q_0 =
\begin{cases}
\mathbb{H}^{m-1} \perp \langle 1, -c \rangle, & \text{if } n = 2m, \\
\mathbb{H}^{m} \perp \langle -c \rangle, & \text{if } n = 2m+1,
\end{cases}
\]
with $c \in k^*$ (see \cite[Table p.55]{Tits1966} for $B_n$ and see \cite[Table p.56]{Tits1966} for $D_n$).

\medskip

We recall the following elementary lemma on algebraic groups :

\begin{lemma} \label{lemma_cover}
    Let $k$ be a field, not necessarily of characteristic zero, let $G$ be a semisimple quasi-split algebraic $k$-group and $\widetilde{G}$ its simply connected cover.
    Then $H^1(k,\widetilde{G}) \longrightarrow H^1(k,G)$ is injective.
\end{lemma}
\begin{proof}  The kernel $\mu$ of the map $\widetilde{G} \longrightarrow G$ is central.  
Since $\widetilde{G}$ is quasi-split, semisimple, and simply connected, by \cite[XXIV Prop.~3.13]{SGA3} there exists a quasi-split torus $T \subset \widetilde{G}$, that is, of the form $R_{A/k}(\mathbb{G}_m)$ for some $k$-étale algebra $A$. Since $C_G(T)=T$ and $\mu$ is central, it is contained in $T$ . The kernel of the map
\[
H^1(k, \widetilde{G}) \longrightarrow H^1(k, G)
\] 
is the image of  
\[
H^1(k, \mu) \longrightarrow H^1(k, \widetilde{G}).
\]  
This map factors through $H^1(k, T) = 0$ and we get the desired result.
\end{proof}

\begin{proposition}
Question \ref{conjecture} admits a positive answer for semisimple simply connected quasi-split groups of type $B_n$ ($n \geqslant 2$) and $D_n$ ($n \geqslant 4$) defined over a field $k$ of characteristic zero, provided that $k$ has finite cohomological dimension or contains a square root of $-1$.
\end{proposition}

\begin{proof} If $k$ contains a square root of $-1$, we may assume that $A$, $U$, $G$, and $X$ are defined over a subfield of $k$ finitely generated over $\mathbb{Q}(\sqrt{-1})$. We may therefore reduce to the case where $k$ has finite cohomological dimension.

Let $q_0$ be a hyperbolic form (if $G$ is of type $B_n$) or one of the quadratic forms described above (if $G$ is of type $D_n$).  
Let $U \subset A$ be an open subset as in Question~\ref{conjecture}, and let
$\alpha \in H^1(U, \mathrm{Spin}(q_0))$
denotes the class of a $\mathrm{Spin}(q_0)$-torsor over $U$ that is trivial at the origin.
Under the composite
\[
H^1(U, \mathrm{Spin}(q_0)) \longrightarrow H^1(U, \mathrm{SO}(q_0) \longrightarrow W(k(A)),
\]
the image of $\alpha$ is the classe of a quadratic form $q_{\alpha}$ such that the difference $q_{\alpha} - q_0$ lies in the third power of the fundamental ideal (see \cite[§7.1]{Serre_1993_1994}); that is,
\[
q_{\alpha} - q_0 \in I^3(k(A)) \subset W(k(A)).
\]
Moreover, since both $q_{\alpha}$ and $q_0$ are defined over the open $U$ of $A$ and $\mathrm{codim}_A(A\setminus U) \geqslant 2$, their difference lies in the Witt group of the local ring at every codimension-one point of $A$, that is,
\[
q_{\alpha} - q_0 \in W(\mathcal{O}_{A,x}) \quad \text{for all } x \in A^{(1)}.
\]
By Lemma~\ref{exacte_sequence_Witt_group}, it follows that
\[
q_{\alpha} - q_0 \in W_{\mathrm{nr}}(k(A)) \cap I^3(k(A)) = I^3_{\mathrm{nr}}(k(A)) \subset I_{\mathrm{nr}}(k(A)).
\]
Then, by Proposition~\ref{canceling_unramified_witt_group}, and by the same specialisation arguments involving diagram \ref{Witt_diagram} in the proof of Proposition~\ref{conjecture_A_C}, there exists an étale isogeny $f$ of $A$ such that
\[
f^*(q_{\alpha} - q_0) = 0.
\]
By Lemma \ref{lemma_cover}, we see that the kernel of the natural map  
\[
i \colon H^1\!\left(k(A), \mathrm{Spin}(q_0)\right) \;\longrightarrow\; H^1\!\left(k(A), \mathrm{SO}(q_0)\right)
\]  
is trivial.
It follows that $f^*(\alpha)_{|k(A)} = 0$ in $H^1(k(A), \mathrm{Spin}(q_0))$.
This completes the proof.
\end{proof}

\textit{$\bullet$ The cases $G_2$, $F_4$, $E_6$, $E_7$, $^3D_4$, $^6D_4$:}

\medskip

For exceptional groups other than $E_8$, we consider the \emph{Rost invariant} as defined in \cite[Appendix~B.6]{Esnault1998}:
\[ r \colon H^1(X, G) \longrightarrow \mathbb{H}^4_{\mathrm{\acute{e}t}}(X, \Gamma_X(2)),
\]
where $X$ is a smooth variety over a field $F$, and $G$ is a semisimple, simply connected group and $\Gamma_X(2)$ denotes Lichtenbaum’s complex (see \cite{Lichtenbaum1984}, \cite{Lichtenbaum1987}, \cite{Lichtenbaum1990}). More precisely, the authors construct an element
$\gamma_G \in \mathbb{H}^4_{\mathrm{\acute{e}t}}(BG, \Gamma_{BG}(2))$,
where $BG$ denotes the simplicial classifying space of $G$ defined in \cite[§~4.2]{Esnault1998}. For every $G$-torsor $E$ over $X$, one has $r([E]) = \varphi_E^*(\gamma_G)$, where $[E]$ denotes the class of $E$ in $H^1(X,G)$ and 
$\varphi_E \colon X \longrightarrow BG$ is the classifying morphism of $E$, uniquely determined up to homotopy. In particular $r$ is functorial with respect to arbitrary morphisms between smooth varieties. 
We have the following exact sequence:

\begin{proposition}{\cite[Th.~1.1]{Kahn1996}}\label{exact_sequence_gamma}
    Let $X$ be a smooth scheme over a field $F$. We have a short exact sequence
    \[
    0 \longrightarrow \CH^2(X) \longrightarrow \mathbb{H}_{\mathrm{\acute{e}t}}^4(X, \Gamma_X(2)) \longrightarrow H^0_{\mathrm{Zar}}(X, \mathcal{H}^3(\mathbb{Q}/\mathbb{Z}(2))) \longrightarrow 0
    \]
    which is functorial in $X$.
\end{proposition}

\begin{remark}\label{Rost_field}
    In the case where $X = \operatorname{Spec}(F)$ for a field $F$, Proposition~\ref{exact_sequence_gamma} yields an isomorphism
    \[
    \mathbb{H}_{\mathrm{\acute{e}t}}^4(F, \Gamma_{\Spec F}(2)) \cong H^3(F, \mathbb{Q}/\mathbb{Z}(2)),
    \]
    and the Rost invariant coincides with the one originally defined by Rost over fields (see \cite[§7.3]{Serre_1993_1994}).
\end{remark}

\bigskip

The Rost invariant allows us to answer the question positively in the case of exceptional groups $G_2$, $F_4$, $E_6$, $E_7$, $^3D_4$ and $^6D_4$ since the map
\[
H^1(F, G) \longrightarrow H^3(F, \mathbb{Q}/\mathbb{Z}(2))
\]
has trivial kernel when $G$ is an almost simple simply connected quasi-split group of exceptional type distinct from $E_8$ (see \cite[§8 and §9]{Serre_1993_1994} for $G_2$ and $F_4$, \cite[Main Theorem 0.1]{Garibaldi2001} or \cite[Th. 6.1]{Chernousov2003} for $E_6$ and $E_7$ and \cite{Garibaldi2010} for $^3D_4$ and $^6D_4$).


%

\begin{proposition}
    Question \ref{conjecture} admits a positive answer for semisimple simply connected quasi-split groups of types $G_2$, $F_4$, $E_6$, $E_7$, ${}^3D_4$, and ${}^6D_4$ over a field $k$ of characteristic $0$ which contains a square root of $-1$.
\end{proposition}

\begin{proof}
Since $k$ contains a square root of $-1$, we may assume that $A$, $U$, $G$, and $X$ are defined over a subfield of $k$ finitely generated over $\mathbb{Q}(\sqrt{-1})$. We may therefore reduce to the case where $k$ has finite cohomological dimension.

Let $G$ be a semisimple, simply connected, quasi-split group of type 
$G_2$, $F_4$, $E_6$, $E_7$, $^3D_4$, or $^6D_4$. Let $U \subset A$ be an open subset as in Question~\ref{conjecture}, and let 
$\alpha \in H^1(U, G)$ be trivial at the origin. 
There exists an integer $n \in \{2,3,6,12\}$ such that, for every field extension $F/k$, the image of the Rost invariant of $G$ is contained in $H^3(F,\mu_n^{\otimes 2})$ (see \cite[Rem.~10.8, p.~133]{Garibaldi2003}). Moreover, for every field extension $F/k$,
\[
H^3(F, \mu_n^{\otimes 2}) \hookrightarrow H^3(F, \mathbb{Q}/\mathbb{Z}(2))
\]
(see \cite[§7.3]{Serre_1993_1994}), so the Rost invariant takes values in $H^3(F,\mu_n^{\otimes 2})$.  
Since $\mathrm{codim}_A(A\setminus U) \geqslant 2$, by \cite[Thm.~11.1]{Garibaldi2003}, applied to the completion of $k(A)$ at each discrete valuation corresponding to a codimension $1$ point of $A$, it follows that 
$r(\alpha_{|k(A)}) \in H^3_{\mathrm{nr}}(A,\mu_n^{\otimes 2})$.
Since $-1$ is a square in $k$, we have $\mu_n^{\otimes 2} \cong \mathbb{Z}/n$ for all $n \in \{2,3,6,12\}$. It follows from Proposition~\ref{cancelling_unramified_coho} that there exists an étale isogeny 
$f \colon A \to A$ such that $r(f^*\alpha_{|k(A)}) = (f^*r(\alpha))_{|k(A)}$
is a constant element of $H^3_{\mathrm{nr}}(A,\mu_n^{\otimes 2})$.  
Since $\alpha$ is trivial at the origin, and since the Rost invariant and the exact sequence~\ref{exact_sequence_gamma} are functorial, it follows that $r(f^*\alpha)_{|k(A)}=0$.
Since the Rost invariant has trivial kernel for almost simple, simply connected groups of exceptional types distinct from $E_8$ 
(see \cite{Serre_1993_1994, Garibaldi2001, Chernousov2003, Garibaldi2010}), 
it follows that the image of $f^*\alpha$ in $H^1(k(A), G)$ is trivial. 
This concludes the proof.
\end{proof}

From the discussion above, we deduce the following proposition:

\begin{proposition} \label{conjecture_sssc} 
Question~\ref{conjecture} admits a positive answer for semisimple simply connected quasi-split groups whose simple factors are not of type $E_8$, over any field $k$ of characteristic zero that contains a square root of $-1$.
\end{proposition}

\begin{proof}
    The result follows from the decomposition theorem for semisimple simply connected groups (see \cite[(26.8)]{Knus_et_al_1998}) and the preceding discussion.
\end{proof}

\begin{remark}
    In dimension $2$, as mentioned above, the theorem is known without any assumption on the characteristic of the field or on the root datum of the semisimple simply connected group (see \cite[Theorem~1.4]{Jong2011}).
\end{remark}

\subsubsection*{The case of semisimple groups defined over \texorpdfstring{$k$}{}}

We now turn to the slightly more general case of torsors under quasi-split semisimple groups defined over $k$.

\begin{proposition} \label{conjecture_semisimple_k}
   If $k$ is a field of characteristic zero containing a square root of~$-1$, then Question~\ref{conjecture} admits a positive answer for quasi-split semisimple groups over~$k$ whose root data have no factors of type~$E_8$.

If $\dim A = 2$ and $k$ is algebraically closed, then the question admits a positive answer for all semisimple groups over~$k$.

\end{proposition}

\begin{proof}
Let $U \subset A$ be an open subset as in Question~\ref{conjecture} and let $G$ be a semisimple quasi-split group whose root datum contains no factor of type~$E_8$ if $\dim A > 2$ or if $k$ is not algebraically closed.

\noindent
Let $\alpha \in H^1(U, G)$ be the class of a $G$-torsor that is trivial at the origin.
Let $\widetilde{G}$ denotes the semisimple simply connected cover of $G$.
Then we have a short exact sequence
\[
1 \longrightarrow \mu \longrightarrow \widetilde{G} \longrightarrow G \longrightarrow 1,
\]
where $\mu$ is a finite diagonalizable group scheme.
This exact sequence induces, for every open subset $V \subset A$, an exact sequence in cohomology
\[
H^1(V, \widetilde{G})
\longrightarrow
H^1(V, G)
\longrightarrow
H^2(V, \mu).
\]
By Corollary~\ref{cancelling_etale_coho}, since $\alpha$ is trivial at the origin, there exists an étale isogeny
$f \colon A \to A$ such that the image of the class $f^*\alpha$ in
$H^2(f^{-1}(U), \mu)$ is trivial .
Therefore, after pulling back along the étale isogeny $f$, replacing
$U$ by $f^{-1}(U)$ (which still contains $0$ and all points of codimension~$1$,
since $f$ is a finite flat morphism of group schemes), and replacing $\alpha$
by $f^*\alpha$, we may assume that $\alpha$ lifts to a class $\gamma \in H^1(U, \widetilde{G})$.
Let $\gamma_0 := [0]^*\gamma \in H^1(k,\widetilde{G})$.
Let $E_0$ be a $\widetilde{G}$-torsor over $\Spec k$ with class $\gamma_0$.
By \cite[§~2.6.1]{Giraud1971}, there exists an isomorphism
\[
\tau_{E_0} \colon H^1(U, \widetilde{G}^{E_0}) \xlongrightarrow{\sim} H^1(U, \widetilde{G})
\]
such that $\tau_{E_0}^{-1}(\gamma)$ is trivial at the origin.
Moreover, since the image of $\gamma_0$ in $H^1(k, G)$ is $[0]^*\alpha$, which is trivial by hypothesis, there exists a central isogeny $\widetilde{G}^{E_0} \to G$ whose kernel is of multiplicative type; in particular, $\widetilde{G}^{E_0}$ is quasi-split.
Thus, by Proposition~\ref{conjecture_sssc} or Proposition~\ref{thm:serre_conj_II} in the case where the field is algebraically closed and $\dim A = 2$, after pulling back along a suitable isogeny and possibly replacing $U$, we may assume that $\gamma_{|k(A)} = s^{*}(\gamma_0)$,
where $s : \Spec k(A) \to \Spec k$ is the morphism induced by the structure morphism $A \to \Spec k$.Finally, since the image of $\gamma_0$ in $H^1(k, G)$ is trivial, $\alpha_{|k(A)} \in H^1(k(A), G)$ is also trivial, as it is the image of $\gamma_{|k(A)}$. This concludes the proof.

\end{proof}

\subsubsection*{Case of tori defined over $A$}

Before treating the general case of a reductive group, we first consider the case of a torus, which only requires results of section \ref{section_etale_coho_finite_coeff} to be proved.

\begin{proposition} \label{conjecture_torus} Let $k$ be a field of characteristic zero. 
Question~\ref{conjecture} admits a positive answer for any torus $T$ defined over $A$.
\end{proposition}

\begin{proof} Let $U \subset A$ be an open subset as in Question~\ref{conjecture}.  
Let $T$ be a torus defined over $U$, and let $\alpha \in H^1(U, T)$ be a class that is trivial at the origin.  
By Corollary~\ref{torus_constant},after pulling back by an isogeny of $A$ we may assume that $T$ is defined over $k$.  
Then there exists a nonzero integer $n$ such that, for every field extension $k'/k$,  
\[
n\, H^1(k', T) = 0.
\]
Consider the short exact sequence
\[
0 \longrightarrow {}_n T \longrightarrow T \xlongrightarrow{[n]} T \longrightarrow 0,
\]
which gives rise to the exact sequence in cohomology:
\[
H^1(V, T) \xlongrightarrow{[n]} H^1(V, T) \longrightarrow H^2(V, {}_n T)
\]
for any open subset $V\subset A$.
Since ${}_n T$ is a finite multiplicative group, Corollary~\ref{cancelling_etale_coho} implies that there exists an étale isogeny $f$ of $A$ such that the image of $\alpha$ in $H^2(f^{-1}(U), {}_n T)$ becomes constant, and hence trivial, as $\alpha$ is trivial at the origin by hypothesis.  
Therefore, after pulling back along this étale isogeny $f$, the class $f^*\alpha$ is a multiple of $n$ in $H^1(f^{-1}(U), T)$, and is thus generically trivial by the choice of $n$.
\end{proof}

\subsubsection*{Case of reductive groups defined over $A$}

We are now ready to give a positive answer to Question~\ref{conjecture} in the case where $G$ is a reductive group over $U \subset A$, with $U$ an open subset that contains all points of codimension $1$, and $0 \in U$, such that the fiber of $G$ at the origin is a quasi-split reductive group whose root datum contains no factors of type $E_8$.

Since $G$ is connected, it has a well-defined root datum, and hence there exists a unique Chevalley group $\widetilde{G}$ over $\mathbb{Z}$ corresponding to this root datum (see \cite[XXV, §1.2]{SGA3}).  
We then have the following property:

\begin{proposition}  \label{étale_form} Let $k$ be a field not necessarily of characteristic zero. Let $S$ be a connected scheme over $k$ and $s\in S(k)$. let $G$ be a (connected) reductive group scheme over $S$ and $G_s$ it fiber at $s$. Then the group $G$ is an étale form of $G_0 \times_k S$.
\end{proposition}

\begin{proof}
Since $G$ is connected, it has a well-defined root datum, and hence there exists a unique Chevalley group $\widetilde{G}$ over $\mathbb{Z}$ corresponding to this root datum (see \cite[XXV, §1.2]{SGA3}).  
Both $G$ and $G_0 \times_{\Spec k} U$ are forms of $\widetilde{G} \times U$.  
This follows from the fact that $G$ is locally split for the étale topology (see \cite[XXV.~§2.3]{SGA3}), together with the equivalence between split reductive groups and root data (see \cite[XXV.~§1.1]{SGA3}).  
From this, we deduce that $G$ is an étale form of $G_0 \times_{\Spec k} S$.
\end{proof}

We can now prove the following proposition:

\begin{proposition}\label{isomorphism_generic} Let $k$ be a field of characteristic zero that contains a square root of $-1$.  
Let $A$ be an abelian variety over $k$, and let $U \subset A$ be an open subset that contains all points of codimension $1$, with $0 \in U$.  
Let $G$ be a reductive group over $U$ such that its fiber at the origin, denoted $G_0$, is quasi-split and whose root datum has no factors of type~$E_8$ whenever $\dim A > 2$ or $k$ is not algebraically closed.  
Then there exists a class
\[
\gamma \in H^1(f^{-1}(U), (G_0)_{\mathrm{ad}})
\]
 such that
 \begin{enumerate}
     \item it maps to the class of $f^*G$ in $H^1(f^{-1}(U),\mathrm{Aut}(G_0))$,
     \item  $\gamma$ is generically trivial
     \item $\gamma$ is trivial at the origin.
 \end{enumerate} 
\end{proposition}

\begin{proof}
Since $G$ is a form of $G_0 \times U$, it corresponds to an element $\alpha \in H^1(U, \mathrm{Aut}(G_0))$ which is trivial at the origin.  
We have a short exact sequence
\[
0 \longrightarrow (G_0)_{\mathrm{ad}} \longrightarrow \mathrm{Aut}(G_0) \longrightarrow \mathrm{Out}(G_0) \longrightarrow 0,
\]
which induces, for every open subset $V \subset A$, the exact sequence of pointed sets
\[
H^1(V, (G_0)_{\mathrm{ad}}) \longrightarrow H^1(V, \mathrm{Aut}(G_0)) \longrightarrow H^1(V, \mathrm{Out}(G_0)).
\]
The group $\mathrm{Out}(G_0)$ is a twisted constant sheaf over $k$. Hence, by Corollary~\ref{cancelling_etale_coho_non_commutative}, since $\alpha$ is trivial at the origin, after pulling back along an étale isogeny and possibly replacing $U$, we may assume that the image of $\alpha$ in $H^1(U, \mathrm{Out}(G_0))$ is trivial. In other words, $\alpha$ is inner, that is, $G$ is an inner form of $G_0 \times U$.  
Let $\gamma \in H^1(U, (G_0)_{\mathrm{ad}})$ be a class mapping to $\alpha$.  
Let $\gamma_0 := 0^*\gamma$, and let $E_0$ be a $(G_0)_{\mathrm{ad}}$-torsor over $\Spec k$ whose class in $H^1(k, (G_0)_{\mathrm{ad}})$ is $\gamma_0$.  
Let $\Gamma_0$ denotes the pullback of $\gamma_0$ to $H^1(U, (G_0)_{\mathrm{ad}})$.  
Since the image of $\gamma_0$ in $H^1(k, \mathrm{Aut}(G_0))$ is trivial, the same holds for the image of $\Gamma_0$ in $H^1(U, \mathrm{Aut}(G_0))$, so $\Gamma_0$ arises from an element $c \in \mathrm{Out}(G_0)$.  
By \cite[Prop.~3.3.3]{Giraud1971}, there is an action of $c$ on the fiber of the map
\[
H^1(U, (G_0)_{\mathrm{ad}}) \longrightarrow H^1(U, \mathrm{Aut}(G_0)),
\]
and this action sends $\gamma$ to a class that is trivial at the origin.  
Hence, we may choose $\gamma$ so that it is trivial at the origin.  
Since $G_0$ is defined over $k$, the group $(G_0)_{\mathrm{ad}} \times_k U$ is also defined over $k$ and is semisimple.  
Therefore, by Proposition~\ref{conjecture_semisimple_k}, after a further pullback via an étale isogeny $f$ of $A$, the associated class $\gamma$ becomes generically trivial.  
It follows that, after this pullback, the group $f^*G$ is generically isomorphic to $G_0 \times \Spec k(A)$.
\end{proof}

Let us recall the following point of the theory of torsors (see \cite[\S 2.6.3]{Giraud1971}):

\begin{lemma}\label{basic_fact_torsor}
    Let $X$ be a scheme, and let $H$ be a group scheme over $X$. Suppose that $H'$ is an inner form of $H$, defined by an $H_{\mathrm{ad}}$-torsor $P$ over $X$, and assume that $P$ lifts to an $H$-torsor $\widetilde{P}$. Then there is a natural bijection of sets
\[
    \tau_{\widetilde{P}} \colon H^1(X, H') \xlongrightarrow{\sim} H^1(X, H),
\]
which sends the class of the trivial $H'$-torsor to the class $[\widetilde{P}]$ of $\widetilde{P}$. This bijection depends only on the class of $P'$ in $H^1(X, H)$, and it is an isomorphism of pointed sets when $P'$ is trivial.
\end{lemma}

\begin{theorem} \label{Theorem_torsor} Let  $k$ be a field of characteristic zero that contains a square root of $-1$.
Question~\ref{conjecture} admits a positive answer for (connected) reductive groups $G$ defined over $A$, whose root datum contains no factors of type $E_8$ and whose fiber at the origin $G_0$ is quasi-split.

If $\dim A = 2$ and $k$ is algebraically closed, Question~\ref{conjecture} admits a positive answer for all (connected) reductive groups $G$ defined over $A$.
\end{theorem} 

\begin{proof} 
Let $G$ be a reductive group defined over $U$, whose root datum contains no factors of type~$E_8$ whenever $\dim A > 2$ or when $k$ is not algebraically closed.

\textit{First step:} Suppose first that $G$ has finite center, i.e., that $G$ is semisimple. Then $G_0$ also has finite center, and so does $G_0 \times U$. By Proposition~\ref{isomorphism_generic}, after pulling back by an étale isogeny, we may assume that $G$ is an inner form of $G_0 \times U$.  
Let $\gamma \in H^1(U, (G_0)_{\mathrm{ad}})$ denote a class associated with this inner form.  
By \emph{loc.~cit.}, we may assume that $\gamma$ is trivial at the origin and generically trivial.

The following sequence
\[
H^1(V, G_0) \longrightarrow H^1(V, (G_0)_{\mathrm{ad}}) \longrightarrow H^2(V, Z(G_0))
\]
is exact for every open subset $V \subset A$.  
By Proposition~\ref{cancelling_etale_coho}, after pulling back $G$ along an étale isogeny $f$ of $A$, replacing $U$ by $f^{-1}(U)$ and $\gamma$ by $f^*\gamma$, we may assume that the class $\gamma \in H^1(U, (G_0)_{\mathrm{ad}})$ lifts to a class $\tilde{\gamma} \in H^1(U, G_0)$. Let $\tilde{\gamma}_0 = 0^* \tilde{\gamma} \in H^1(k, G_0)$.  
Let $\tilde{E}_0$ be a $G_0$-torsor whose corresponding class in $H^1(k, G_0)$ is $\tilde{\gamma}_0$.  
Set $G'_0 := (G_0)^{\tilde{E}_0}$.  
Since the image of $\tilde{\gamma}_0$ in $H^1(k, (G_0)_{\mathrm{ad}})$ is trivial, $G'_0 \cong G_0$, and hence $G'_0$ is still quasi-split.  
Moreover, $G$ is an inner form of $G'_0 \times U$, as both are inner forms of $G_0 \times U$.  
Let $\tilde{F}_0 := \tilde{E}_0 \wedge^{G'_0} (G'_0)_{\mathrm{ad}}$ and define $\widetilde{\gamma} := \tau^{-1}_{\tilde{F}_0}(\gamma)$; then $\widetilde{\gamma} \in H^1(U,(G'_0)_{\mathrm{ad}})$ is a class associated with the inner form $G$ of $G'_0 \times U$.  
By setting $\tilde{\delta} := \tau^{-1}_{\tilde{E}_0}(\tilde{\gamma})$, we obtain a lift of this class in $H^1(U, \widetilde{G}_0)$ which is trivial at the origin.  
By Proposition~\ref{conjecture_semisimple_k}, we may assume, after possibly pulling back along an étale isogeny $f$ of $A$ and replacing $U$, that the class
\[
\tilde{\delta} \in H^1(U, \widetilde{G}_0)
\]
is generically trivial, i.e., trivial over the generic point $\Spec k(A)$.  
Let $P'$ be a $\widetilde{G}_0$-torsor representing $\tilde{\delta}$, and denote by $P'_g$ its restriction to $\Spec k(A)$.  
Fix an isomorphism
\[
G \;\overset{\phi}{\cong}\; {}^{P'}(G'_0 \times A).
\]
Then, by Lemma~\ref{basic_fact_torsor}, we obtain the commutative diagram
\begin{equation} \label{diagram_commutative_form}
\begin{tikzcd}[sep=huge]
H^1(U, G)  
  \arrow[r, "\sim","\phi"'] \arrow[d] 
& H^1(U, {}^{P'}(G'_0 \times A))  
  \arrow[r, "\sim","\tau_{P'}"'] \arrow[d]  
& H^1(U, G'_0) \arrow[d] \\
H^1(k(A), G)  
  \arrow[r, "\sim","\phi_{k(A)}"']  
& H^1(k(A), {}^{P'}(G'_0 \times A))  
  \arrow[r, "\sim","\tau_{P'_g}"'] 
& H^1(k(A), G'_0)
\end{tikzcd}
\end{equation}
in which all horizontal arrows are bijections.  
Moreover, the bottom horizontal arrow is an isomorphism of pointed sets, since $\tilde{\delta}$ is generically trivial (see Lemma~\ref{basic_fact_torsor}).

Now let $\alpha \in H^1(U, G)$ be a class trivial at the origin.  
As $\tilde{\delta}$ is trivial at the origin, $(\tau_{P'} \circ \phi)(\alpha)$ is also trivial at the origin.  
By Proposition~\ref{conjecture_semisimple_k}, there exists an étale isogeny $f$ of $A$ such that the pullback
\[
f^*\bigl(\tau_{P'} \circ \phi(\alpha)\bigr) \in H^1(f^{-1}(U), G'_0)
\]
is generically trivial.  
By commutativity of diagram~\eqref{diagram_commutative_form}, it follows that $f^*\alpha$ is generically trivial.  
This completes the first step.

\medskip

\textit{Second step:} Consider now the general case, without assuming that $G$ has finite center. Let $\alpha \in H^1(U, G)$ be trivial at the origin. By \cite[XXII Th.~6.2.1]{SGA3}, we have the exact sequence of group schemes:
\begin{equation}\label{exact_sequence_DG}
    0 \longrightarrow DG \longrightarrow G \longrightarrow \operatorname{corad}(G) \longrightarrow 0,
\end{equation}
and the map $G \longrightarrow \operatorname{corad}(G)$ restricts to a smooth isogeny $\operatorname{rad}(G) \longrightarrow \operatorname{corad}(G)$.  
By Corollary~\ref{torus_constant}, after pulling back along a suitable isogeny of $A$ and replacing $U$, we may assume that $\operatorname{corad}(G)$ and $\operatorname{rad}(G)$ are defined over $k$. Let $\mu$ denotes the finite multiplicative group given by the kernel of the map $\operatorname{rad}(G) \longrightarrow \operatorname{corad}(G)$.  
Since $U$ is quasi-compact, there exists a non-zero integer $n$ such that $n\mu = 0$.  
As $\operatorname{rad}(G)$ is central, pulling back \eqref{exact_sequence_DG} along 
\[
\operatorname{corad}(G) \xlongrightarrow{n} \operatorname{corad}(G)
\] 
yields the following commutative diagram:
\[
\begin{tikzcd}[sep=huge]
1 \arrow[r] & DG \arrow[r] & G \arrow[r] & \operatorname{corad}(G) \arrow[r] & 1 \\
1 \arrow[r] & DG \arrow[r] \arrow[u] & DG \times \operatorname{corad}(G) \arrow[r] \arrow[u] & \operatorname{corad}(G) \arrow[r] \arrow[u, "n"'] & 1
\end{tikzcd}
\]
Thus, we obtain the exact sequence
\[
0 \longrightarrow {}_n\mathrm{corad}(G) \longrightarrow DG \times \operatorname{corad}(G) \longrightarrow G \longrightarrow 0,
\]
where ${}_n\mathrm{corad}(G)$ is the kernel of $\operatorname{corad}(G) \xlongrightarrow{n} \operatorname{corad}(G)$.
This induces the corresponding cohomology exact sequence:
\[
    H^1(U, DG \times \operatorname{corad}(G)) \longrightarrow H^1(U, G) \longrightarrow H^2(U, {}_n\mathrm{corad}(G)).
\]
Since ${}_n\mathrm{corad}(G)$ is finite and of multiplicative type, and since $\alpha$ is trivial at the origin, it follows from Proposition~\ref{cancelling_etale_coho} that, after pulling back along an étale isogeny of $A$, we may assume that $\alpha$ is the image of a class $\widetilde{\alpha} \in H^1(U, DG \times \operatorname{corad}(G))$.  
As $H^1(U, DG \times \operatorname{corad}(G))= H^1(U,DG) \times H^1(U,\operatorname{corad}(G))$, we write $\widetilde{\alpha}=(\widetilde{\alpha}_1,\widetilde{\alpha}_2)$. Moreover, since $\alpha$ is trivial at the origin and $\operatorname{corad}(G)$ is defined over $k$, 
\[
    0^* \widetilde{\alpha} \in H^1\big(k, (DG \times \operatorname{corad}(G)) \times_A \Spec k\big)
\] 
comes from a class 
\[
    \widetilde{\alpha}'_0 \in H^1(k, {}_n\mathrm{corad}(G)).
\] 

Let $P''_0$ be a torsor whose class is $\widetilde{\alpha}'_0$, and set 
\[
    P'' := P''_0 \times A.
\]
Define
\[
    P' := P''  \wedge^{{}_n\mathrm{corad}(G)} (DG \times \operatorname{corad}(G)).
\] 
Then $P'$ is isomorphic to the product $P'_1 \times_A P'_2$, where $P'_1$ is a $DG$-torsor and $P'_2$ is a $\operatorname{corad}(G)$-torsor over $A$.  

The image of $\tau_{P'}^{-1}(\widetilde{\alpha})$ in 
\[
H^1(U, {}^{P'}(DG \times \mathrm{corad}(G))) 
\;\cong\; H^1(U, {}^{P'_1} DG) \times H^1(U, {}^{P'_2} \operatorname{corad}(G))
\] 
is trivial at the origin.  
Since the image of $P'$ in $H^1(A,G)$ is trivial, there exists a multiplicative isogeny ${}^{P'}(DG \times \mathrm{corad}(G)) \longrightarrow G$,  
so that the fiber at the origin of ${}^{P'}(DG \times \mathrm{corad}(G))$ is quasi-split, and hence the fiber at the origin of ${}^{P'_1} DG$ is quasi-split.  
Therefore, by the first step and Proposition~\ref{conjecture_torus}, we may pull back $\tau_{P'}^{-1}(\widetilde{\alpha})$ along a suitable étale isogeny of $A$ so that it becomes generically trivial.  
It follows that, after pulling back $\widetilde{\alpha}$ along an étale isogeny $f$ of $A$, the class $f^* \widetilde{\alpha}$ is generically equal to the class associated with the torsor $P' \times \Spec k(A)$.  

Since 
\[
    P' \times \Spec k(A) = (P'' \times \Spec k(A)) \wedge^{{}_n\mathrm{corad}(G) \times \Spec k(A)} ((DG \times \operatorname{corad}(G)) \times_A \Spec k(A)), 
\] 
the image of $(f^* \widetilde{\alpha})_{k(A)} \in H^1(k(A), DG \times \operatorname{corad}(G))$ in $H^1(k(A), G)$, i.e., $(f^* \alpha)_{k(A)}$, is trivial.  
Therefore, the result follows.
\end{proof}

\subsection{The question for Borel varieties} \label{TheconjectureforBorelvarieties}

We now study the case of smooth projective families of homogeneous spaces of connected reductive groups.

\begin{definition}
Let $X \longrightarrow S$ be a morphism of schemes. We say that $X$ is a \emph{Borel $S$-scheme} if:
\begin{enumerate}
  \item[\textup{(i)}] $X$ is proper and smooth over $S$;
  \item[\textup{(ii)}] for every $s \in S$, the geometric fiber $X_s$ is a Borel variety.
\end{enumerate}
\end{definition}

Let us recall that to any (connected) reductive group scheme $G$ over $S$, one can associate the scheme $P(G)$ of parabolic subgroups (see \cite[XXVI Thm.~3.3]{SGA3}).  
There is a morphism
\[
P(G) \longrightarrow \mathrm{Of}(\mathrm{Dyn}(G))
\]
which associates to each (smooth) parabolic subgroup its Dynkin type.  
Here, $\mathrm{Of}(\mathrm{Dyn}(G))$ is the finite étale scheme over $S$ of open and closed subschemes of $\mathrm{Dyn}(G)$ defined in \cite[XXVI §3.1]{SGA3}.

Given a type $\tau$, we define the scheme \( P_\tau(G) \) as the fiber product:
\[\begin{tikzcd}[sep=huge]
P_\tau(G) \arrow[r] \arrow[d]  & P(G) \arrow[d] \\
S \arrow[r,"\tau"] & \mathrm{Of}(\mathrm{Dyn}(G)).
\end{tikzcd}\]

Borel varieties exhibit a remarkably rigid structure.

\begin{proposition}\cite[Prop.~4]{Demazure1977} \label{Demazure_Borel}
Let $X$ be a Borel $S$-scheme. Then  $\underline{\mathrm{Aut}}_S(X)$ is a semisimple $S$-group scheme, and $X$ is a homogeneous space under of $G$, the connected component of $\underline{\mathrm{Aut}}_S(X)$. In particular, $X$ is $G$-isomorphic to a scheme of the form $P_\tau(G)$ for some type $\tau$.
\end{proposition}

\begin{remark}
\begin{enumerate}
    \item[\textup{(i)}]  Proposition~\ref{Demazure_Borel} shows that such a family is always projective (see \cite[XXVI~Th.~3.3(ii)]{SGA3}).
    \item[\textup{(ii)}]  As the center of $G$ acts trivially on $X$, the group $G$ is adjoint.
    
    \item[\textup{(iii)}]  If $S$ is connected, the group $G$ has constant type, given by a Dynkin diagram that we call the \emph{type} of $X$.
\end{enumerate}
\end{remark}

This rigidity allows us to answer Question~\ref{conjecture_originale} positively, away from the $E_8$ case, and without significant extra work in a slightly more general setting: namely, we allow the field $k$ to have cohomological dimension at most $1$.

\begin{theorem}\label{Theorem_homogeneous} Let $k$ be a field of characteristic zero, of cohomological dimension at most~$1$ and containing a square root of~$-1$. Let $A$ be an abelian variety over a $k$ and let $U \subset A$ be an open subset that contains all points of codimension $1$, with $0 \in U$.
Let $X \longrightarrow U$ be a smooth, projective family of homogeneous spaces under connected reductive algebraic groups.  
In each of the following cases:
\begin{enumerate}
    \item[\textup{(i)}] $\dim A = 2$ and the field~$k$ is algebraically closed;
    
    \item[\textup{(ii)}] the type of $X$ does not contain factors of type~$E_8$,
\end{enumerate}
there exists an étale isogeny $f$ of~$A$ such that the pullback
\[
f^*X \longrightarrow f^{-1}(U)
\]
admits a rational section.
\end{theorem}

\begin{proof}
   By Proposition~\ref{Demazure_Borel}, $X$ is of the form $P_\tau(G)$ for some connected semisimple group $G$ defined over $U$.  
Since $\mathrm{cd}(k) \leqslant 1$, the fiber of $G$ at the origin is quasi-split by Steinberg's theorem (see \cite[Cor.~5.2.6]{Gille2019}).  
Hence, by Proposition~\ref{isomorphism_generic}, we may pull back $G$ along an étale isogeny so that
\[
    G \times_A \Spec k(A) \cong G_0 \times \Spec k(A),
\]
where $G_0$ denotes the fiber of $G$ at the origin.  
Since reductive groups over a field of cohomological dimension at most $1$ admit parabolic subgroups of any given type (again by Steinberg's theorem), we conclude that 
\[
    P_\tau(G) \times_A \Spec k(A) \cong  P_\tau(G_0) \times_A \Spec k(A)
\] 
has a $k(A)$-point.
\end{proof}

\begin{remark}\label{remark_U}
    It is immediate that, if we allow pullbacks by isogenies of $A$ which may not be group morphisms, Theorem~\ref{Theorem_torsor} and Theorem \ref{Theorem_homogeneous} still holds under the hypothesis that $U \subset A$ is an open subset such that $\mathrm{cd}(A \setminus U) \geqslant 2$, equipped with a $k$-point $u_0$, and that the fiber of $G$ at $u_0$ is quasi-split.
\end{remark}

\subsection{Some applications}\label{Generalisation}

Theorems~\ref{Theorem_torsor} and~\ref{Theorem_homogeneous} can be generalized to the case where the torsor or the homogeneous space is defined not over an abelian variety but over a Borel $A$-scheme.  
\smallskip

Before extending these results, we recall a well-known fact about torsors.

\begin{lemma}\label{acceptable_group} Let $k$ be a field of characteristic zero. 
Let $H$ be a locally constant group scheme over $k$, and let $G$ be a reductive algebraic group over $k$. 
Let $G'$ be an algebraic group over $k$ fitting into an exact sequence
\[
1 \longrightarrow G \longrightarrow G' \longrightarrow H \longrightarrow 1.
\]
Then:
\begin{enumerate}
    \item[\textup{(i)}] There is a canonical isomorphism
    \[
        H^1(k,G') \;\cong\; H^1(\mathbb{A}^1_k,G').
    \]

   \item[\textup{(ii)}] For every integer $N \ge 1$ and every $k$-point, 
\[
\ker\!\Bigl(H^1(\mathbb{A}^N_k,G') \longrightarrow H^1(k,G')\Bigr)
\;\subset\;
H^1_{\mathrm{Zar}}(\mathbb{A}^N_k,G').
\]
\end{enumerate}
\end{lemma}

\begin{proof}
(i) We first prove the statement for a locally constant group $H$ over $k$.  
By \cite[Lem.~2.13(2)]{Gille24}, every $H$-torsor over $\mathbb{A}^1_{k_s}$ is isotrivial and hence trivial, since $\mathbb{A}^1_{k_s}$ is simply connected (see \cite[Th.~1.1]{KambayashiSrinivas1983}).  
Thus, by \cite[§2.2]{Gille2015}, every element of $H^1(\mathbb{A}^1_k,H)$ comes from an element of 
\[
H^1(\mathrm{Gal}(k_s/k),H(\mathbb{A}^1_{k_s})) 
   \;=\; H^1(\mathrm{Gal}(k_s/k),H(k_s)) 
   \;=\; H^1(k,H),
\]
so the result holds for locally constant $H$. 

\noindent Let $G'$ be as in the statement, and let $\alpha \in H^1(\mathbb{A}^1_k, G')$.
We want to show that $\alpha$ is constant.  
Let $E$ be the $G'$-torsor over $\mathbb{A}^1_k$ corresponding to the class $\alpha$.
Consider the pullback of $E$ along the zero section $0 : \Spec k \to \mathbb{A}^1_k$, $E' := E \times_{\mathbb{A}^1_k} \Spec k$.
There is a natural morphism of algebraic groups over $k$
\[
(G')^{E'} \longrightarrow H^{E'},
\]
whose kernel is a reductive algebraic group $G''$ over $k$.
The image of $\alpha$ in
\[
H^1(\mathbb{A}^1_k, (G')^{E'})
\]
under the map $\tau_{E'}^{-1}$ comes from an element of
\(H^1(\mathbb{A}^1_k, G'')\).
Indeed, $\tau_{E'}^{-1}(\alpha)$ is trivial at the origin, and every
$H^{E'}$-torsor is constant by the preceding step.
By \cite{Raghunathan1984}, any $G''$-torsor over $\mathbb{A}^1_k$ that is
trivial at a point is constant. Hence the class in
\(H^1(\mathbb{A}^1_k, G'')\) coming from $\tau_{E'}^{-1}(\alpha)$ is constant.
Since $E'$ comes from a constant torsor, and the image of $\alpha$ in
\(H^1(\mathbb{A}^1_k, (G')^{E'})\) is constant, it follows that
the original class \(\alpha \in H^1(\mathbb{A}^1_k, G')\) is constant.

\noindent (ii) This is an consequence of \cite[Th.~A]{Raghunathan1978} together with~(i).
\end{proof}

We now generalise our result.

\begin{proposition}
Let $k$ be an algebraically closed field of characteristic zero and let $A$ be an abelian variety defined over~$k$.  
Let $G$ be a reductive group over an irreducible Borel $A$-scheme~$X$.  
Assume that one of the following conditions holds:
\begin{enumerate}
    \item[\textup{(i)}] $\dim A = 2$;
    \item[\textup{(ii)}] both the type of $G$ and the type of $X$ have no factors of type~$E_8$.
\end{enumerate}  
Then there exists a reductive group $G_0$ defined over~$k$ such that $G$ is a form of $G_0 \times A$, and there exists a finite étale covering
\[ f \colon X' \longrightarrow X \]
of~$A$ such that, after pulling back $G$ via~$f$, the group $G \times_X X'$ is generically isomorphic to $G_0 \times \Spec k(X')$.  
Moreover, the corresponding class $\gamma$ can be chosen to be trivial at the origin, and $X'$ can be taken to be a Borel $A$-scheme.
\end{proposition}

\begin{proof}
Let $\phi \colon X \longrightarrow A$ denote the structure morphism.  
By Theorem~\ref{Theorem_homogeneous}, after possibly pulling back along an étale isogeny of~$A$, we may assume that $\phi$ admits a section generically.  
Moreover, a closer inspection of the proof shows that we can choose this étale isogeny so that, after the pullback, the geometric generic fiber of~$\phi$, namely
\[
X_{k(A)} := X \times_A \Spec k(A),
\]
is of the form $G'_{k(A)}/P'_{k(A)}$, where $G'$ is a reductive group defined over~$k$ and $P' \subset G'$ is a parabolic subgroup.  
The generic section can then be taken as the base change of a section $\Spec k \longrightarrow G'/P'$ of the morphism $G'/P' \longrightarrow \Spec k$.  
By the Bruhat decomposition (see \cite[XXVI~4.3.2~(iv')]{SGA3}), we may assume that our section
\[
s \colon \Spec k(A) \longrightarrow X_{k(A)}
\]
factors through an affine open subset
\[
\mathbb{A}^N_{k(A)} \subset X_{k(A)}.
\]
Since $\phi$ is proper, this section is defined over an open subset $U \subset A$ containing all points of codimension~$1$.  
We continue to denote this section by
\[
s \colon U \longrightarrow X.
\]
Let $u_0$ be a $k$-point of~$U$, and set $x_0 := s(u_0)$.  
Let $G_0$ denote the fiber of~$G$ at~$x_0$.  
As $X$ is irreducible, $G$ has constant type, equal to that of~$G_0$.  
Hence, by Proposition~\ref{étale_form}, the group $G$ is an étale form of $G_0 \times X$.
Let $\alpha \in H^1(X, \mathrm{Aut}(G_0))$ be the class corresponding to~$G$.  
We can pull back this class to $U$ via the section $s$.  
Then, by Proposition~\ref{isomorphism_generic}, Remark~\ref{remark_U}, and our choice of $s$, after pulling back along an étale isogeny of $A$ (which is not necessarily a group morphism), we may assume that there exists a section
\[
s \colon \Spec k(A) \longrightarrow X
\]
such that the pullback of $\alpha$ by $s$ is trivial in $H^1(k(A), \mathrm{Aut}(G_0))$, and that the restriction of this section factors through an affine open subset
\[
\mathbb{A}^N_{k(A)} \subset X_{k(A)}.
\]
The kernel of the map
\[
H^1\!\bigl(\mathbb{A}^N_{k(A)}, \mathrm{Aut}(G_0)\bigr)
\longrightarrow
H^1\!\bigl(k(A), \mathrm{Aut}(G_0)\bigr)
\]
is contained in 
\(H^1_{\mathrm{Zar}}\!\bigl(\mathbb{A}^N_{k(A)}, \mathrm{Aut}(G_0)\bigr)\)
by Proposition~\ref{acceptable_group}(ii).  
It follows that $\alpha$ is generically locally trivial. Since we have only performed finite étale base changes on~$A$,  
we may take the étale covering $X' \longrightarrow X$ that makes $G$ generically isomorphic to $G_0 \times \Spec k(X')$ to be an irreducible Borel $A$-scheme.
\end{proof}

This property enables us to prove, exactly as in Theorem~\ref{Theorem_homogeneous}, the following result:

\begin{theorem} \label{Theorem_rational_section} Let $A$ be an abelian variety over an algebraically closed field $k$ of characteristic zero, and let $X$ be an irreducible Borel $A$-scheme.  
Let $Y \longrightarrow X$ be a smooth, projective family of homogeneous spaces under connected reductive algebraic groups.  
Assume one of the following conditions holds:
\begin{enumerate}
  \item[\textup{(i)}] $\dim A = 2$;
  \item[\textup{(ii)}] the type of $X$ and the types of $Y$ contain no factors of type~$E_8$.
\end{enumerate}
Then there exists an étale cover $X' \longrightarrow X$, with $X'$ an irreducible Borel $A$-scheme, such that the pullback family $Y \times_X X' \longrightarrow X'$ admits a rational section.
\end{theorem}

We now generalize Theorem~\ref{Theorem_torsor} to the case of $G$-torsors over Borel $A$-varieties, restricting ourselves to reductive groups $G$ defined over~$k$.

\begin{theorem}\label{theorem_Borel}
 Let $A$ be an abelian variety over an algebraically closed field $k$ of characteristic zero, and let $X$ be an irreducible Borel $A$-scheme.  
Let $G$ be a reductive group over~$A$, and let $Y \longrightarrow X$ be a $G$-torsor.  
In each of the following cases:
\begin{enumerate}
  \item[\textup{(i)}] $\dim A = 2$;
  \item[\textup{(ii)}] the type of $X$ and the type of $G$ contain no factors of type~$E_8$,
\end{enumerate}
there exists a finite étale cover $X' \to X$ with $X'$ an irreducible Borel $A$-scheme such that $Y \times_X X'$ is Zariski-locally trivial.
\end{theorem}

\begin{proof}
Let $\phi \colon X \longrightarrow A$ denote the structure morphism.  
By Theorem~\ref{Theorem_homogeneous}, after possibly pulling back along an étale isogeny of~$A$, we may assume that $\phi$ admits a section generically.  
Moreover, a closer inspection of the proof shows that we can choose this étale isogeny so that, after the pullback, the geometric generic fiber of~$\phi$, namely
\[
X_{k(A)} := X \times_A \Spec k(A),
\]
is of the form $G'_{k(A)}/P'_{k(A)}$, where $G'$ is a reductive group defined over~$k$ and $P' \subset G'$ is a parabolic subgroup. %
The generic section can then be taken as the base change of a section $\Spec k \longrightarrow G'/P'$ of the morphism $G'/P' \longrightarrow \Spec k$.
By the Bruhat decomposition (see \cite[XXVI~4.3.2~(iv')]{SGA3}), we may assume that our section
\[
s:\Spec k(A) \longrightarrow X_{k(A)}
\]
factors through an affine open subset
\[
\mathbb{A}^N_{k(A)} \subset X_{k(A)}.
\]
Since $\phi$ is proper, this section extends to an open subset $U \subset A$ containing all points of codimension~$1$.  
We continue to denote this section by
\[
s \colon U \longrightarrow X.
\]

\noindent Let $Y_U := Y \times_X U$.  
By Theorem~\ref{Theorem_torsor} and Remark~\ref{remark_U}, there exists an étale isogeny $f \colon A \longrightarrow A$
not necessarily a morphism of groups, such that the pullback of $Y_U$ along~$f$ is Zariski-locally trivial over $f^{-1}(U)$.  
Thus, after pulling back by an étale isogeny we may assume that $Y_U$ admits a rational section.
Hence, we obtain a section
\[
\Spec k(A) \xlongrightarrow{s_{k(A)}} X_{k(A)} \longrightarrow \Spec k(A),
\]
such that $s_{k(A)}^{*}Y$ is trivial and such that $s_{k(A)}$ factors through an affine open subset
\[
\mathbb{A}^N_{k(A)} \subset X_{k(A)}.
\]
 As
\[
\ker\!\left(H^{1}\bigl(\mathbb{A}^N_{k(A)}, G\bigr)
  \longrightarrow H^{1}\bigl(k(A), G\bigr)\right)
  \subset H^{1}_{\mathrm{Zar}}\bigl(\mathbb{A}^N_{k(A)}, G\bigr)
\]
by Proposition~\ref{acceptable_group}(ii). We deduce that the $G$-torsor~$Y$ is rationally trivial.  
Finally, since all finite étale coverings considered arise from base change by isogenies of~$A$ (not necessarily morphisms of groups), we may take the cover $X' \longrightarrow X$ appearing in the statement to be an irreducible Borel $A$-scheme.
\end{proof}
We now state some consequences for smooth projective varieties with nef tangent bundles.
We first recall the relevant definitions and one of the main structure theorems for such varieties.

\begin{definition} 
    Let $X$ be a projective scheme over a field $k$ of characteristic zero.
    \begin{enumerate}
        \item[\textup{(i)}] A line bundle $\mathcal{L}$ on $X$ is said to be \emph{nef} if for any curve $C \subset X$ one has
        \[
        \mathcal{L} \cdot C \geqslant 0.
        \]
        \item[\textup{(ii)}] A vector bundle $E$ on $X$ is said to be \emph{nef} if the line bundle $\mathcal{O}_{\mathbb{P}(E)}(1)$ is nef on $\mathbb{P}(E)$.
    \end{enumerate}
\end{definition}

There is a structure theorem for varieties with nef tangent bundle:

\begin{theorem}[\cite{DemaillyPeternellSchneider1994}, Main Theorem]\label{Demailly_Theorem}
Let $k$ be an algebraically closed field of characteristic zero, and let $X$ be a smooth, projective, connected $k$-scheme with nef tangent bundle.
Let
$\pi \colon \widetilde{X} \longrightarrow X$
be an étale cover of maximal irregularity $q(\widetilde{X}) := h^{1}\bigl(\widetilde{X}, \mathcal{O}_{\widetilde{X}}\bigr)$.
Then the scheme $\widetilde{X}$ admits a fibration given by its Albanese map
\[
\widetilde{X} \longrightarrow \operatorname{Alb}(\widetilde{X}),
\]
where $\operatorname{Alb}(\widetilde{X})$ denotes its Albanese variety, which is an abelian variety.
Moreover, the fibres of this fibration are Fano varieties with nef tangent bundle.
\end{theorem}

We now state the following conjecture.

\begin{conjecture}[\cite{CampanaPeternell1991}, Conjecture~11.1]
Over an algebraically closed field of characteristic zero, a Fano variety with nef tangent bundle is a projective homogeneous variety of a linear algebraic group, i.e.\ of the form $G/P$ for $G$ a connected linear algebraic group and $P$ a parabolic subgroup.
\end{conjecture}

\begin{remark}\label{remark_C_P}
This conjecture was proved by Campana and Peternell in dimension less than~$3$ (see \cite[Th.~p.~169]{CampanaPeternell1991}) and by J.-M.~Hwang \cite{Hwang2006} in dimension~$4$. 
More generally, see \cite{MunozOcchettaSolaWatanabeWisniewski2015} for a survey on the conjecture.
\end{remark}

Under the assumption of the \emph{Campana--Peternell conjecture}, the finite étale cover $\widetilde{X}$ of a projective variety with nef tangent bundle and maximal irregularity is, by Theorem~\ref{Demailly_Theorem}, an $\operatorname{Alb}(\widetilde{X})$-Borel variety.
Theorem~\ref{Theorem_homogeneous} then has consequences for the fibration of Theorem~\ref{Demailly_Theorem}:

\begin{corollary} Let $k$ be a field of characteristic zero.
Let $X$ be a smooth projective $k$-scheme with nef tangent bundle.
Let 
\[
\widetilde{X} \longrightarrow \operatorname{Alb}(\widetilde{X})
\] 
be the fibration given by Theorem~\ref{Demailly_Theorem}.
Assume that the \emph{Campana--Peternell Conjecture} holds for the fibres of 
$\widetilde{X} \longrightarrow \operatorname{Alb}(\widetilde{X})$, for instance when $\dim X \leqslant 4$ (see Remark~\ref{remark_C_P}), so that $\widetilde{X}$ becomes an $\operatorname{Alb}(\widetilde{X})$-Borel variety.
If one of the following condition is satisfied
\begin{enumerate}
    \item[\textup{(i)}] the type of $\widetilde{X}$ contains no factors of type~$E_8$;
    \item[\textup{(ii)}] $\dim \operatorname{Alb}(\widetilde{X}) = 2$,
\end{enumerate}
then there exists an étale isogeny 
$f \colon \operatorname{Alb}(\widetilde{X}) \longrightarrow \operatorname{Alb}(\widetilde{X})$
such that the pullback family
$f^*\widetilde{X} \longrightarrow \operatorname{Alb}(\widetilde{X})$
admits a rational section.
\end{corollary}

Furthermore, Theorem~\ref{theorem_Borel} yields consequences for $G$-torsors over smooth complex varieties with nef tangent bundles.

\begin{corollary} Let $k$ be field of characteristic zero.
Let $X$ be a smooth projective $k$-scheme with nef tangent bundle.  
Let $G$ be a reductive group over $k$, and let $Y \longrightarrow X$ be a $G$-torsor.
Let 
\[
\widetilde{X} \longrightarrow \operatorname{Alb}(\widetilde{X})
\] 
be the fibration given by Theorem~\ref{Demailly_Theorem}.
Assume that the \emph{Campana--Peternell Conjecture} holds for the fibres of 
$\widetilde{X} \longrightarrow \operatorname{Alb}(\widetilde{X})$, for instance when $\dim X \leqslant 4$ (see Remark~\ref{remark_C_P}), so that $\widetilde{X}$ becomes an $\operatorname{Alb}(\widetilde{X})$-Borel variety.
If one of the following condition is satisfied
\begin{enumerate}
 \item[\textup{(i)}] the type of $\widetilde{X}$ and the type of $G$ contain no factors of type~$E_8$;
    \item[\textup{(ii)}] $\dim \operatorname{Alb}(\widetilde{X}) = 2$,  
\end{enumerate}
then there exists a finite étale cover $X' \longrightarrow X$, with $X'$ a smooth projective variety with nef tangent bundle, such that the pullback $Y \times_X X'$ is Zariski-locally trivial.
\end{corollary}

\appendix
\section*{\normalfont \textsc{Appendix}}

\renewcommand{\thetheorem}{A.\arabic{theorem}}
\setcounter{theorem}{0}

In this appendix, we recall basic results of intersection theory that are used in the first paragraph to lift the results of Moonen and Polishchuk on the Chow group of the Jacobian of curves to the motives of the Jacobian of curves.

Let $f : X \longrightarrow Y$ be a morphism between smooth projective schemes.
To $f$, we can associate its graph $[\Gamma_f] \in \CH_*(Y \times_k X)$, denoted $[f]_*$, as well as its transpose ${}^t[\Gamma_f] \in \CH_*(X \times_k Y)$ denoted $[f]^*$. We then have the following result.

\begin{proposition} \label{action_graphe}
Let $S$ be a smooth projective scheme over $k$.
Let $\alpha \in \CH_*(X \times S)$ and $\beta \in \CH^*(Y \times S)$. Then:
\begin{enumerate}
   \item[\textup{(i)}]  $\alpha \circ [f]_* = (f_{S})_*(\alpha)$,
   \item[\textup{(ii)}] $\beta \circ [f]^* = (f_{S})^*(\beta)$.
\end{enumerate}
\end{proposition}

\begin{proof}
\begin{enumerate}
   \item[\textup{(i)}] Follows from \cite[Prop.~16.1.1(c)(i)]{Fulton1998}.
   \item[\textup{(ii)}] Follows from \cite[Prop.~16.1.1(c)(ii)]{Fulton1998}, together with \cite[Prop.~16.1.1(b)]{Fulton1998}.
\end{enumerate}
\end{proof}

\begin{remark}
By definition, $[f]^*$ induces a morphism from $h(Y)$ to $h(X)$, while $[f]_*$ induces a morphism from $h(X)^{\vee}$ to $h(Y)^{\vee}$ (see \cite[§4.1]{André2004}).
\end{remark}

\begin{lemma} \label{intersection_theory_lemma}
Let $X$ and $Y$ be two smooth projective schemes over $k$, and let $f : Y \longrightarrow X$ be a morphism of schemes.
Let $S'$ be a smooth projective scheme over $k$, let $c \in \Corr(S, S') = \CH_*(S \times_k S')$, and let $a \in \CH_*(X \times_k S)$ and $b \in \CH_*(Y \times_k S)$. Then:
\begin{enumerate}
    \item[\textup{(i)}] $c \circ (f_S^* a) = f_{S'}^*(c \circ a)$ in $\CH_*(Y \times_k S')$,
    \item[\textup{(ii)}] $c \circ (f_S)_* b = (f_{S'})_*(c \circ b)$ in $\CH_*(X \times_k S')$.
\end{enumerate}
\end{lemma}

\begin{proof}
(i) We compute:
\begin{align*}
c \circ (f_S^* a)
&= c \circ (a \circ [f]_*) 
&& \hspace{-2.54cm} \text{\footnotesize (ii)~\ref{action_graphe}} \\[4pt]
&= (c \circ a) \circ [f]_*
&& \hspace{-2.54cm} \text{\footnotesize \cite[Prop.~16.1.1~(a)]{Fulton1998}} \\[4pt]
&= f_{S'}^*(c \circ a)
&& \hspace{-2.54cm} \text{\footnotesize (ii)~\ref{action_graphe}}
\end{align*}

(ii) Similarly,
\begin{align*}
c \circ ((f_S)_* b)
&= c \circ (b \circ [f]_*) 
&& \hspace{-2.54cm} \text{\footnotesize (i)~\ref{action_graphe}} \\[4pt]
&= (c \circ b) \circ [f]_*
&& \hspace{-2.54cm} \text{\footnotesize \cite[Prop.~16.1.1~(a)]{Fulton1998}} \\[4pt]
&= (f_{S'})_*(c \circ b)
&& \hspace{-2.54cm} \text{\footnotesize (i)~\ref{action_graphe}}
\end{align*}
\end{proof}

\renewcommand*{\bibfont}{\small}
\printbibliography

\bigskip

\textsc{Institut Camille Jordan, Université Claude Bernard Lyon 1, Bâtiment Jean Braconnier 21 Av. Claude Bernard, 69100 Villeurbanne and Institute of Mathematics ”Simion Stoilow” of the Romanian Academy, 21
Calea Grivitei Street, 010702 Bucharest, Romania}

\emph{Email address:} \texttt{bruneaux@math.univ-lyon1.fr}

\end{document}